\documentclass[12pt]{article}

\usepackage{latexsym,amsmath,amscd,amssymb,graphics}
\usepackage{enumerate,framed,comment}
\usepackage{color,url}
\usepackage{verbatim}
\usepackage[margin=2.5cm]{geometry}
\usepackage{amsthm}
\usepackage{epsfig,epstopdf,color,enumerate}
\usepackage{caption}
\usepackage{float}
\usepackage{placeins}
\usepackage{graphicx, subfigure, psfrag} 

\usepackage[colorlinks]{hyperref}
\usepackage{url}  %
\usepackage[percent]{overpic}

\usepackage[T1]{fontenc}
\usepackage[colorlinks]{hyperref}
\usepackage{sidecap}

\usepackage[all]{xy}

\DeclareMathOperator{\Id}{Id}
\DeclareMathOperator{\ad}{ad}

\DeclareMathOperator{\cay}{cay}

\newcommand{\rem}[1]{}
\makeatletter

\@addtoreset{figure}{section}
\def\thefigure{\thesection.\@arabic\c@figure}
\def\fps@figure{h, t}
\@addtoreset{table}{bsection}
\def\thetable{\thesection.\@arabic\c@table}
\def\fps@table{h, t}
\@addtoreset{equation}{section}

\makeatother

\rem{
\textwidth 6.5 truein
\oddsidemargin 0 truein
\evensidemargin .2 truein
\topmargin -.6 truein
\textheight 9.1 in
}

{

\begin{document}

\newtheorem{theorem}{Theorem}[section]
\newtheorem{definition}[theorem]{Definition}
\newtheorem{lemma}[theorem]{Lemma}
\newtheorem{remark}[theorem]{Remark}
\newtheorem{proposition}[theorem]{Proposition}
\newtheorem{corollary}[theorem]{Corollary}
\newtheorem{example}[theorem]{Example}

\def\below#1#2{\mathrel{\mathop{#1}\limits_{#2}}}



\title{
Geometric Integrators for Higher-order\\ Mechanics on Lie Groups}
\author{
Christopher L. Burnett$^{1}$, Darryl D. Holm$^{1}$, David M. Meier$^{1}$
}
{\addtocounter{footnote}{1} 
\footnotetext{Department of Mathematics, Imperial College, London SW7 2AZ, UK. 
\texttt{c.burnett@imperial.ac.uk, d.holm@ic.ac.uk, d.meier09@ic.ac.uk}
}%
%
\date{}

\maketitle

\makeatother
\maketitle


%

\begin{abstract} 
This paper develops a structure-preserving numerical integration scheme for a class of higher-order mechanical systems. The dynamics of these systems are governed by invariant variational principles defined on higher-order tangent bundles of Lie groups. The variational principles admit Lagrangians that depend on acceleration, for example. The symmetry reduction method used in the Hamilton--Pontryagin approach for developing variational integrators of first-order mechanics is extended here to higher order. The paper discusses the general approach and then focuses on the primary example of Riemannian cubics. Higher-order variational integrators are developed both for the discrete-time integration of the initial value problem and for a particular type of trajectory-planning problem. The solution of the discrete trajectory-planning problem for higher-order interpolation among points on the sphere illustrates the approach.
\end{abstract}
\newpage
\tableofcontents


\section{Introduction}\label{Intro-sec}
After introducing the general background for higher-order mechanics and variational integrators, we state our goals and summarize the content of the paper. 
\subsection{General background}
The problem treated in this paper fits into a classical type of problem in control theory called \emph{trajectory planning}, or \emph{interpolation by variational curves}. The task in this type of problem is to find an optimal curve that interpolates through a given set of points (or configurations) lying in a manifold. The configuration manifold is specific to the application. In this paper, we shall study the example of trajectory planning for tracking a rigid body through a prescribed sequence of orientations. The control of a rigid body summons either the group $SO(3)$ of rotations in $\mathbb{R}^3$, or the semidirect-product group $ SE(3) = SO(3)\rtimes \mathbb{R}^3$ of three-dimensional rotations and translations in Euclidean space. Such trajectory-planning problems are relevant in numerous applications, for example, in aeronautics, robotics, computer-aided design, air traffic control, biomechanics and more recently, computational anatomy. 

\paragraph{Higher order mechanics.} Some trajectory-planning applications may require the optimal trajectories to possess a certain degree of smoothness. This requirement summons variational principles that depend on higher-order derivatives of the interpolation path, such as \emph{acceleration} (rate of change of velocity) or \emph{jerk} (rate of change of acceleration) etc. Properties of these higher-order variational principles, including their Hamiltonian formulation and their symmetry reduction, have been studied in \cite{deLRo1985, BedD2005, Gay-BalmazEtAl2011HOLPHP, Gay-BalmazEtAl2010, CodD2011}.
 An interesting example of such a higher-order variational principle was introduced in \cite{GaKa1985} and \cite{NoHePa1989}. This example requires the minimization of the mean-square covariant acceleration and leads to curves called \emph{Riemannian cubic splines}. These curves generalize the familiar cubic splines of Euclidean space to the context of Riemannian manifolds. The mathematical theory of Riemannian cubics and their higher-order generalizations has been developed in a series of papers \cite{CrSL1995, CaSLCr1995, CaSLCr2001, No2004, No2006, splinesanalyse}. Applications to computer graphics, spacecraft control and computational anatomy are discussed in \cite{PaRa1997, ZeKuCr1998, HuBl2004, TrVi2010}, amongst others. We refer to \cite{Popiel2007, MaSLeKr2010} and \cite{Noakes_SphericalSplines2006} for extensive references and historical discussions concerning Riemannian cubics, their higher-order generalizations, and related higher-order interpolation methods.

Many properties of first-order mechanics remain relevant in the higher-order setting. In particular, curves that solve Hamilton's principle for a given non-degenerate Lagrangian in first-order mechanics have various interesting geometric features. For example, the corresponding Hamiltonian and symplectic structure, both defined on the cotangent bundle of the configuration manifold, will be conserved along these solutions. In addition, if the Lagrangian is invariant under the action of a Lie group then the associated \emph{momentum map} is also conserved (we refer to \cite{MaRa03} for a detailed discussion). In fact, much of this structure has been extended to higher order \cite{deLRo1985, Gay-BalmazEtAl2010, Gay-BalmazEtAl2011HOLPHP, CodD2011}. Moreover, in first order mechanics much has also been achieved in bringing these properties to the numerical side, where time is discretised and one aims to approximately solve the variational equations. The numerical schemes oriented towards preserving such properties form the topic of \emph{geometric integration}, which is a relatively new subject for higher-order mechanics \cite{HLW06,CoJidD2011}.
\paragraph{Variational integrators.} 
In designing a geometric integrator, one must face the fact that all three properties (momentum preservation, symplecticity and energy preservation) cannot be achieved simultaneously by methods with a fixed time step \cite{KaMaOr99,Ma_LOM,HLW06}. Our focus is on the symplectic and momentum preserving schemes that follow from discretising Hamilton's variational principle. The resulting variational integrators are well established \cite{Cadzow70, Maeda82, WeMa97} and the discrete analogues of many of the familiar properties implied by Hamilton's principle, such as Noether's theorem, have been established \cite{Logan73}. This progress has made variational integration a convenient approach for numerics associated with mechanics \cite{Lee87}. Despite lacking energy preservation, variational integrators that are symplectic benefit from the favourable energy behaviour and long time reliability properties afforded by symplectic methods for Hamiltonian systems \cite{Rowlands91, BeGi94, Reich99, HLW06}. The background and breadth of variational integration for Lagrangian mechanics on manifolds is reviewed by Marsden and West \cite{MaWe01}.

For our purposes here, the preservation of Lie group structure is also relevant. This issue has been addressed for many types of geometric integration \cite{LeSi94, HLW06, IMKNZ00, AuKrWa93}. Variational integrators for first-order mechanics on Lie groups and the corresponding discrete Euler--Poincar\'e and Lie--Poisson equations have been discussed, for example, in \cite{MaPeSh99a, BoSu99a, LeLeCl07}. The prevalence of Lie groups in applications means these integrators have been implemented in various fields, including computer graphics, integrable systems, solid mechanics and quantum mechanics \cite{KYTKMSD06, MoVe91, OrSt99, JaNo97a}.

An interesting variant in the approach to these Lie group integrators was introduced in \cite{BRMa09}, where discretisation was applied to the first-order Hamilton--Pontryagin (HP) variational principle \cite{Holm_book2_2008, YoMa06}, rather than to the more familiar Hamilton's principle. The HP principle is of particular significance to the present investigation  because it explicitly contains the definitions of the variables of interest for mechanics, position, velocity and momentum, which are all elements of the Pontryagin bundle \cite{BRMa09}. This explicit dependence in the HP principle allows one to take unconstrained variations over the Pontryagin bundle. This property carries over to the discrete HP principle and simplifies certain aspects of discrete mechanics. For example, it removes the need to identify discrete Legendre transforms when discussing symplectic structure. This property greatly facilitates the move to discrete higher-order mechanics.

\subsection{Motivation} \label{Sec-Motivation}
Our primary motivation is to develop symplectic, momentum-preserving integrators for higher-order mechanics on Lie groups. The methods used to develop Hamilton--Pontryagin integrators in \cite{BRMa09} are particularly amenable to extension to higher order due to the explicit nature of the HP variational principle.  We investigate the geometric properties of the resulting numerical algorithms both for the initial value problem and for applications in trajectory planning. 
\begin{figure}[htb]
\begin{center}
\includegraphics[scale=0.229]{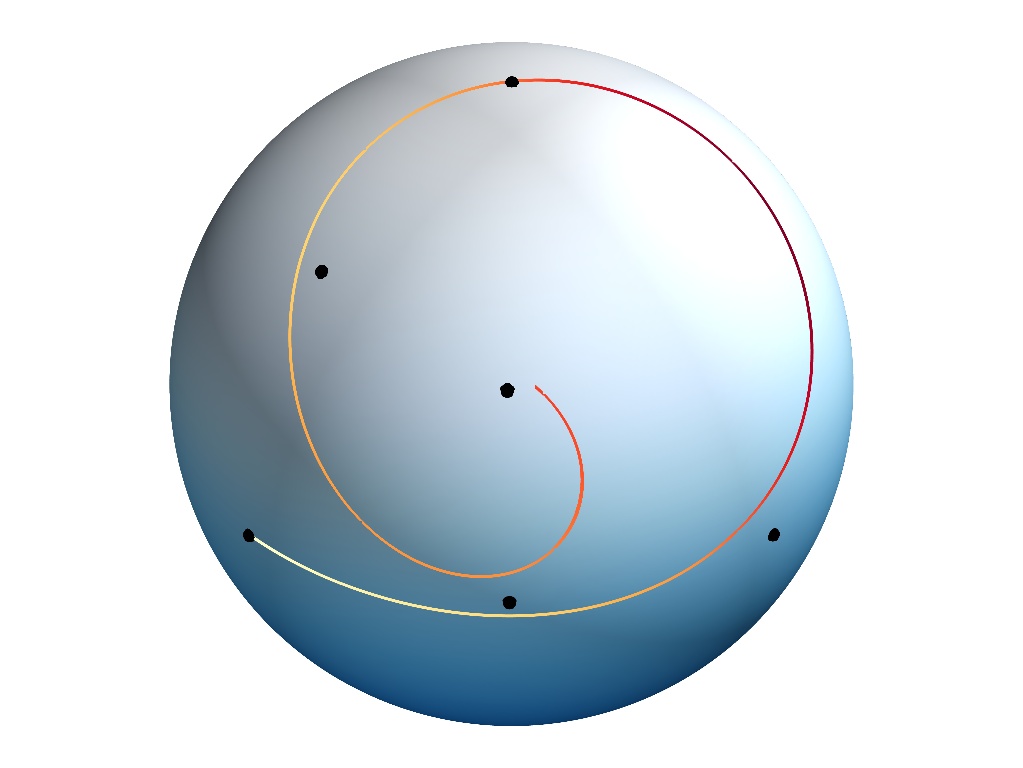}
\caption{\footnotesize Illustration of the trajectory-planning problem on the sphere discussed in Section \ref{applTM-sec}. In this problem, one seeks a curve that passes near the target points (black dots) at prescribed times. The  trajectory shown is generated from a path of rotations acting on Euclidean space. The variational description requires the minimization of a functional that measures both the amount of acceleration (measured in the group of rotations) along the curve and the amount of mismatch between the targets and the curve at the prescribed times. The discussion in Section \ref{applTM-sec} shows that the geometric properties of optimal curves are inherited by the discrete higher-order Hamilton-Pontryagin description of the problem.}
\label{intro_TM_curve}
\end{center}
\end{figure}
\subsection{Main content of the paper} 

The plan of the paper is as follows. 
\vspace{-0.5mm}
\begin{itemize}
\item[] After introducing the basic definitions that we will need, Section \ref{continutime} develops a theorem that characterises Hamilton's principle for second order mechanics on Lie groups (analogous to the first order principle) for Lagrangians with symmetry. We show equivalence with two higher-order versions of variational principles from first-order mechanics: (i) the reduced Hamilton's principle and (ii) the reduced Hamilton--Pontryagin (HP) principle \cite{BRMa09}. In this paper, our main focus is on the latter, in which Lagrange multipliers are used to enforce particular kinematic constraints. We use this principle to derive some previously known results for continuous-time systems, including the NHP equation of \cite{NoHePa1989} for Riemannian cubics. Continuing in Sections \ref{HamStruc} and \ref{Sec-moma_noether} we describe the Hamiltonian theory and symmetry properties underlying higher-order variational problems on Lie groups. In particular, we will discuss how the flow map induced by the HP principle is symplectic and preserves momentum.

\item[] In Section \ref{discretetime} we discuss discrete mechanics and obtain the discrete geometric counterparts of the HP differential equations using a method similar to that in \cite{BRMa09}. The kinematic constraints are discretised using Runge--Kutta methods and their Lie group incarnations, the Runge--Kutta--Munthe-Kaas methods \cite{IMKNZ00}. These are then incorporated by using Lagrange multipliers into a discrete version of the higher-order HP principle. We will find that the flow maps of the resulting algorithms inherit the symplectic momentum-preserving properties of the continuous-time systems in Section \ref{continutime}. Having established the theory, we proceed in Section \ref{AppCubics} with an implementation of it in the particular case of Riemannian cubics on $SO(3)$ with a bi-invariant metric.

\item[]  Section \ref{applTM-sec} applies these discrete methods to a second order trajectory-planning problem in which the interpolating curve evolves by means of a group action.  We first formulate the problem in the framework of the Hamilton--Pontryagin principle.
As in previous work \cite{Gay-BalmazEtAl2010}, we find that the momentum experiences \emph{kicks} related in a simple fashion to the mismatch between the interpolating curve and the target points - a matter which can be presented particularly clearly in the HP framework. We show how the discretisation of the problem respects this geometric feature of solution curves. The section proceeds with a discussion of a numerical algorithm that allows us to calculate curves such as the one shown in Figure \ref{intro_TM_curve} for the trajectory-planning problem on the sphere.
\item[] In Section \ref{Sum-sec} we summarise and present a series of outstanding directions for further research. 
\end{itemize}

\section{Continuous-time mechanics}\label{continutime}
Here we develop the theory of higher-order mechanics in the formalism of the Hamilton-Pontryagin variational principle. After giving some general definitions, we introduce the fundamental dynamical equations. We then proceed to discuss Hamiltonian structure and momentum maps. 
\subsection{Definitions}\label{Sec_Definitions}
We start by introducing the main mathematical objects used in the paper. Let $G$ be a Lie group and $k$ a positive integer, $k\ge0$. The \emph{$k^{th}$-order tangent bundle $\tau_G^{(k)}: T^{(k)}G \rightarrow G$} is defined as a set of equivalence classes of curves, as follows: Two curves $\gamma_i: t \mapsto g_i(t)$, $i = 1, 2$, are \emph{equivalent}, if and only if their time derivatives at $t = 0$ up to order $k$ coincide in any local chart. That is, $g_1^{(l)}(0) = g_2^{(l)}(0)$, for $0\leq l\leq k$. The equivalence class of a curve $\gamma: t \mapsto g(t)$ is denoted by $\left[\gamma\right]_{g(0)}^{(k)}$, or formally as $(g(0), \dot{g}(0), \ldots g^{(k)}(0))$. The set of all equivalence classes of curves based at $g_0$ is written as $T^{(k)}_{g_0}G$.

 The  \emph{$k^{th}$-order Pontryagin bundle} $P^{(2)}G$ is defined as
\cite{CodD2011}
\begin{equation} \nonumber
  P^{(k)}G = T^{(k)}G \times_{T^{(k-1)}G}T^*(T^{(k-1)}G)
\end{equation}
Note that $T^{(0)}G = G$, $T^{(1)}G = TG$ and $P^{(1)}G = TG \times T^*G$. The projection map onto the second factor is denoted by $pr_2: P^{(2)}G \rightarrow T^*(TG)$. For  any fixed real number $a \geq 0$ and smooth manifold $M$, we write $\mathcal{C}(M)$ to denote the set of smooth curves $\gamma~:~ [0, a] \to M$, $t \mapsto  \gamma(t) $. 

The group structure of $G$ affords the following \emph{trivialisations} which make extensive use of the left multiplication maps, $L_g:G\to G$, $h\mapsto gh$. We have $TG \cong G \times \mathfrak{g}$ via $v_g \mapsto  (g, TL_{g^{-1}} v_g) = (g, g^{-1} v_g)$, and $T^*G \cong G \times \mathfrak{g}^*$ via $\alpha_g \mapsto (g, (TL_g)^* \alpha_g) = (g, g^* \alpha_g)$. Here we introduced the notation $ g^{-1} v_g := TL_{g^{-1}} v_g$ and $ g^* \alpha_g: =(TL_g)^* \alpha_g $.  In order to define the trivialisation of $T^{(k)}G$ let $g(t)$ be a representative of $[\gamma]^{(k)}_{g(0)}$ and $\xi(t) = g^{-1}(t) \dot{g}(t)$ be the left-trivialised velocity. The trivialisation map $T^{(k)}G \cong G \times k\mathfrak{g}$ is given by
\begin{equation}
 [\gamma]^{(k)}_{g(0)} \mapsto \left(g(0), \xi(0), \dot{\xi}(0), \ldots, \xi^{(k-1)}(0)\right). \nonumber
\end{equation}
These maps further induce the identification $T^*(T^{(k-1)}G) \cong G \times (k-1) \mathfrak{g} \times k \mathfrak{g}^*$. Accordingly, $P^{(k)}G \cong G \times k\mathfrak{g} \times k\mathfrak{g}^*$. Note that all of these trivialisations make use of the group multiplication from the left. Right-trivialisations can be defined correspondingly using the right multiplication.

From here on we shall consider second order mechanics, $k = 2$, and for later reference we set the notation
\begin{align*}
(g, \xi, u) &\in G \times 2\mathfrak{g} \cong T^{(2)}G\,, \\
(g, \xi, \mu, \nu) &\in G \times \mathfrak{g} \times 2 \mathfrak{g}^* \cong T^*(TG)\,,\\
(g, \xi, u, \mu, \nu) &\in G\times 2\mathfrak{g}\times 2\mathfrak{g}^* \cong P^{(2)}G.
\end{align*}
 We also define the projection $\Pi:  P^{(2)}G \to G$ onto the first factor.
 
\subsection{Variational principles}\label{Sec-Variational_Principles}
A second-order Lagrangian $L: T^{(2)}G \rightarrow \mathbb{R}$ is called $\emph{left-invariant}$ if its left trivialisation $l : G \times 2 \mathfrak{g} \rightarrow \mathbb{R}$ does not depend on the first entry. One then introduces the \emph{reduced} Lagrangian $\ell: 2\mathfrak{g} \rightarrow \mathbb{R}$.
The following theorem, the cornerstone of our approach, can easily be generalised to higher order, $k > 2$.
\emph{
\begin{theorem}\label{Theorem_HOHP}
  Let $L: T^{(2)}G \rightarrow \mathbb{R}$ be a left-invariant Lagrangian with reduced Lagrangian $\ell: 2\mathfrak{g} \rightarrow \mathbb{R}$ and let $g$ be in $\mathcal{C}(G)$. The following are equivalent:
  \begin{enumerate}
  \item Hamilton's principle. The curve $g(t)$ satisfies $\delta S (g)=0$ for $S: \mathcal{C}(G) \rightarrow \mathbb{R}$,
    \begin{equation}
      S = \int_0^a L(g, \dot{g}, \ddot{g})\, dt, \label{Hamiltons_principle}
    \end{equation} 
with respect to arbitrary variations $\delta g\in T_g\mathcal{C}(G)$ which fixed $g(0)$, $g(a)$, $\dot{g}(0)$, $\dot{g}(a)$.
    \item Reduced Hamilton's principle. The curve $v(t):= g^{-1}(t) \dot{g}(t)$ satisfies $\delta \mathrm{s}(v) = 0$ for $\mathrm{s}: \mathcal{C}(\mathfrak{g}) \rightarrow \mathbb{R}$,
      \begin{equation}
        \mathrm{s} = \int_0^a \ell(v, \dot{v})\, dt, \nonumber
      \end{equation}
with respect to variations of the form $\delta v = \dot{\eta} + \operatorname{ad}_{v}\eta$, $\eta\in\mathcal{C}(\mathfrak{g})$ with $\eta(0) = \eta(a) = \dot{\eta}(0) = \dot{\eta}(a) = 0$.
    \item Reduced higher-order Hamilton--Pontryagin (HOHP) principle.
 There exists a curve $\gamma$ in $P^{(2)}G$ with projection $\Pi \circ \gamma = g$ that satisfies $\delta s_{hp}(\gamma) = 0$ for $s_{hp}: \mathcal{C}(P^{(2)}G) \rightarrow \mathbb{R}$,
      \begin{align}
        s_{hp} = \int_0^a \ell(\xi, u) + \left<\mu, g^{-1}\dot{g} - \xi\right>  + \left<\nu, \dot{\xi} - u\right>\, dt,
        \label{HOHPP}
      \end{align}
with respect to variations $\delta \gamma \in T_{\gamma}\mathcal{C}(P^{(2)})$ with $\delta g(0) = \delta g(a) = 0$ and $\delta \xi(0) = \delta \xi(a) = 0$.
  \end{enumerate}
\end{theorem}
}
Equivalence of these variational principles follows from comparing the equations they imply. For illustration, we derive the equations governing solutions to the HOHP principle. For $\gamma\in\mathcal{C}(P^{(2)}G)$ as above denote by $z_\gamma \in T_{\gamma}\mathcal{C}(P^{(2)}G)$ an arbitrary vector over $\gamma$. That is, $z_\gamma: t \mapsto (\delta g(t), \delta \xi(t), \delta u(t), \delta \mu(t), \delta \nu(t))$. We define $\eta = g^{-1} \delta g$ and, upon recalling  $\delta(g^{-1}\dot{g}) = \dot{\eta} + \operatorname{ad}_{g^{-1}\dot{g}}\eta$, we compute
\begin{align}
\begin{split}
\delta s_{hp}=ds_{hp}(z_\gamma) &= \left[ \left< \mu, \eta \right>+\left< \nu, \delta\xi \right>\right]_0^a + \int_0^a \left<\ad^*_{g^{-1}\dot{g}}\mu-\dot{\mu},\,\eta\right>  +\left< \frac{\delta \ell}{\delta \xi}-\mu-\dot{\nu},\,\delta\xi\right> \\
 &\qquad +\left<\frac{\delta \ell}{\delta u}-\nu,\delta u\right> +\left<\delta\mu,g^{-1}\dot{g}-\xi\right>+\left<\delta\nu,\dot{\xi}-u\right>\, dt \,. \label{RHOHP}
\end{split}
\end{align}
The end-point terms will be of interest later, but for now we employ the HOHP principle with $\eta(0) = \eta(a) = \delta \xi(0) = \delta \xi(a) = 0$. We obtain the HOHP equations of motion,
\begin{align}
\dot\mu=&\ad^*_\xi\mu,  \label{RHP_dynamics}\\
\mu=\frac{\delta \ell}{\delta \xi}-\dot{\nu},& \ \nu= \frac{\delta \ell}{\delta u}, \label{Ostograd}\\
\dot{g}=g\xi,& \ \ \dot{\xi}=u, \label{RHP_reconstruction}
\end{align}
which agree with the ones derived, in \cite{Gay-BalmazEtAl2010}, from the reduced Hamilton's principle. Equation \eqref{RHP_dynamics} is usually called the second order Euler--Poincar\'e equations. Equations \eqref{Ostograd} define the Ostrogradsky momenta $\mu$ and $\nu$. 

\paragraph{Riemannian cubics.} The paper focuses on Riemannian cubics. These variational curves were introduced in \cite{GaKa1985} and \cite{NoHePa1989} in the general context of Riemannian manifolds. Here we consider the special case of cubics on Lie groups. 

Let $G$ be a Lie group with a Riemannian metric $d$ and associated norm functions $\left\|.\right\|_g: T_gG \rightarrow \mathbb{R}$. We define the isomorphism $\flat: \mathfrak{g} \rightarrow \mathfrak{g}^*$ by $\left<\eta^\flat, \xi\right> = d(\eta, \xi)$, for all $\xi \in \mathfrak{g}$. Denote its inverse by $\sharp = \flat^{-1}$. Consider the Lagrangian $L: T^{(2)}G \rightarrow \mathbb{R}$,
\begin{equation}\label{Cubic_Lagrangian}
  L(g, \dot{g}, \ddot{g}) = \frac{1}{2} \left\|\frac{D}{Dt}\dot{g}\right\|_g^2,
\end{equation}
where $D/Dt$ is the covariant derivative along curves, with respect to the Levi-Civita connection for the metric $d$. If the metric is left-invariant then $L$ is left-invariant with reduced Lagrangian $\ell: 2\mathfrak{g} \rightarrow \mathbb{R}$,
\begin{equation}\label{Reduced_cubic_Lagrangian}
  \ell(\xi, u) = \frac{1}{2} \left\|u - \operatorname{ad}_\xi \xi^\dagger\right\|_{\mathfrak{g}}^2,
\end{equation}
where $\operatorname{ad}^\dagger$ is the metric adjoint, that is, 
\begin{equation}\label{ad-dagger-def}
\operatorname{ad}^\dagger_\xi \eta := (\operatorname{ad}^*_\xi \eta^\flat)^\sharp
\,,
\end{equation}
for all $\xi, \eta \in \mathfrak{g}$. A \emph{Riemannian cubic} on $G$ is a curve $g(t)$ that is associated with a solution of the HOHP equations \eqref{RHP_dynamics}--\eqref{RHP_reconstruction} for $L$ defined in \eqref{Cubic_Lagrangian}. If the metric is bi-invariant (that is, invariant under both left and right multiplication), then $\operatorname{ad}^\dagger = - \operatorname{ad}$, so that the reduced Lagrangian becomes $\ell(\xi, u) = \frac{1}{2} \left\|u\right\|_{\mathfrak{g}}^2$. In this case the HOHP equations \eqref{RHP_dynamics}--\eqref{RHP_reconstruction} reduce to
\begin{align}
 \dot{\mu}=\ad_{\xi}^* \mu,\quad  \mu=-\dot{u}^{\flat},\quad \nu=u^{\flat},\quad \dot g = g\xi,\quad   \dot{\xi} = u.
\label{HOHPcubic}
\end{align}
These can be rewritten as $\dddot{\xi} = - \operatorname{ad}_\xi \ddot{\xi}$, together with the reconstruction relation $\dot{g} = g\xi$. The solution curves $\xi$ are called \emph{Lie quadratics} in \cite{No2003}. 

For $G = SO(3)$ we obtain the NHP equation, first derived in \cite{NoHePa1989},
\begin{equation}
  \dddot{\xi} = \ddot{\xi} \times \xi. \label{NHP}
\end{equation}
Here we used the identification of the Lie algebra 
 $\mathfrak{so}(3)$ with $\mathbb{R}^3$ together with the cross product, by means of the isomorphism
\begin{align}
\widehat{}:(\mathbb{R}^3,\times)\to\mathfrak{so}(3)\,, \qquad \xi=(x_1,x_2,x_3)\mapsto \widehat{\xi}=\left(\begin{array}{ccc}0 & -x_3 & x_2 \\x_3 & 0 & -x_1 \\-x_2 & x_1 & 0\end{array}\right). \nonumber
\end{align}
\subsection{Hamiltonian structure} \label{HamStruc}
We now demonstrate the symplectic properties of the higher-order dynamics by making explicit use of the HOHP variational principle and \eqref{RHOHP}. This proof will illuminate a similar approach to be used in Section \ref{discretetime} for the corresponding statement in discrete-time mechanics. The results can also be viewed in the context of Hamilton's canonical equations on $T^*(TG)$. 

The second equation in \eqref{Ostograd}  is the constraint relation $\nu = \frac{\delta \ell}{\delta u}$, and we are led to define the submanifold $W_c \subset P^{(2)}G$ that respects this relation. From now on, we assume that the reduced Lagrangian $\ell$ is \emph{hyper-regular}, that is,  $\frac{\delta^2\ell}{\delta u^2}$ is non-degenerate. Therefore, the map $\phi:= \left.pr_2\right|_{W_c}: W_c \rightarrow T^*(TG)$ is a diffeomorphism, with $pr_2$ defined in Section \ref{Sec_Definitions}. We remark that the Lagrangian $\ell = \frac{1}{2} \left\|u\right\|_{\mathfrak{g}}^2$ for Riemannian cubics on manifolds with bi-invariant metrics satisfies hyper-regularity since $\frac{\delta^2\ell}{\delta u^2} = \mbox{Id}$. Here the constraint submanifold $W_c$ is determined by $\nu = u^\flat$. Hyper-regularity is also satisfied when the metric only has one-sided invariance.

 \paragraph{Symplecticity.} Assume that equations \eqref{RHP_dynamics} and \eqref{RHP_reconstruction} induce  \emph{flow maps} $F_t: T^*(TG) \rightarrow T^*(TG)$, for all $t \in [0, a]$. These are defined as follows. For an arbitrary $x \in T^*(TG)$ let $\gamma$ be the solution to the HOHP principle with initial condition $\gamma(0) = \phi^{-1} (x)$. Then by definition $F_t(x) = \phi(\gamma(t))$.

 Next we show that the flow maps preserve the canonical symplectic form on $T^*(TG)  \cong G \times \mathfrak{g} \times 2\mathfrak{g}^*$.  Let $y = (g, \xi, \mu, \nu)$ be in $T^*(TG)$, and let $v_i = (\delta g_i, \delta \xi_i, \delta \mu_i, \delta \nu_i)$, $i = 1, 2$, be arbitrary vectors in $T_y(T^*(TG))$. Define $\eta_i = g_i^{-1}\delta g_i$. The canonical one-form $\theta(y)$ can be calculated as
\begin{equation}\label{Canonical_one_form}
  \theta(y)(v_1) = \left<\mu, \eta_1\right> + \left<\nu, \delta \xi_1\right>,
\end{equation}
and its exterior derivative, the canonical symplectic form, as
\begin{align}
  \omega(y)(v_1, v_2) = -\left<\delta \mu_1, \eta_2\right> - \left<\delta \nu_1, \delta \xi_2\right> + \left<\delta \mu_2, \eta_1\right> + \left<\delta \nu_2, \delta \xi_1\right> + \left<\mu, \left[\eta_1, \eta_2\right]\right>.\label{symp_form}
  \end{align}
The space of solution curves to the HOHP principle can be identified with $T^*(TG)$ by mapping a solution $\gamma$ to $\phi(\gamma(0))$. Correspondingly, the restriction of $s_{hp}$ to the space of solution curves can be identified with a map $s: T^*(TG) \rightarrow \mathbb{R}$. Namely, if $\gamma$ is the solution with initial condition $\gamma(0) = \phi^{-1}(x)$ for $x \in T^*(TG)$, then $s(x) = s_{hp}(\gamma)$. Choose an arbitrary $\delta x \in T_xT^*(TG)$, and set  $\delta \gamma(t)= (\phi^{-1} \circ F_t)_* (\delta x)= (\delta g, \delta \xi, \delta u, \delta \mu, \delta \nu)$, the corresponding variation through solution curves. We compute from  \eqref{RHOHP} and \eqref{Canonical_one_form}
\begin{align}
  ds( \delta x) = \left[\left<\mu, \eta\right> + \left<\nu, \delta \xi\right>\right]^a_0= \left[(F_a)^*\theta - \theta\right] (\delta x), \label{One-form_pull_back}
\end{align}
with $\eta$ as before. The exterior derivative of the \emph{generating function} relation $ds =  (F_a)^*\theta - \theta$ shows that the symplectic form $\omega$ is preserved by $F_a$, and therefore by $F_t$ for all $t \in [0, a]$.

\paragraph{Hamilton's canonical equations.} We close this section with a brief discussion of Hamilton's canonical equations, details of which can be found in \cite{CodD2011} and \cite{Gay-BalmazEtAl2010}. For a hyper-regular Lagrangian $\ell$ the HOHP equations \eqref{RHP_dynamics} and \eqref{RHP_reconstruction} are equivalent to Hamilton's canonical equations
\begin{equation}\label{Hamiltons_canonical_equations}
  i_{X} \omega = dH,
\end{equation}
where $X \in \mathfrak{X}(T^*(TG))$ is the Hamiltonian flow vector field and $H: T^*(TG) \rightarrow \mathbb{R}$ is the Hamiltonian function given by
\begin{equation}\label{Hamiltonian}
  H(x) = \left<\mu, \xi\right> + \left<\nu, u(x)\right> - \ell(\xi, u(x)),
\end{equation}
where $\frac{\delta \ell}{\delta u}( \xi, u(x)) = \nu$. As an example, for cubics
\begin{align}
H(x) = \frac{1}{2}\left\|\nu^\sharp\right\|^2_{\mathfrak{g}}+\left<\mu, \xi\right> .
\label{CubicHam}
\end{align}
As a consequence of \eqref{Hamiltons_canonical_equations}, the Hamiltonian vector field $X$ preserves the canonical symplectic form $\omega$ on $T^*(TG)$. This is an alternative demonstration that $F_t$ is symplectic. 

\subsection{Momentum maps and Noether's theorem}\label{Sec-moma_noether}
In this section we provide a number of necessary definitions concerning group actions and momentum maps before discussing the momentum conservation property of the HOHP flow map arising from the group symmetry of the Lagrangian. 
\paragraph{General definitions.} Consider a left action of a Lie group $H$ on a smooth manifold $Q$,
\begin{equation}
\Phi:H \times Q \rightarrow Q\,, \quad (h, q) \mapsto \Phi_h(q). \nonumber
\end{equation}
 The action, $\Phi:H\times Q\to Q$ induces an infinitesimal generator for each $\zeta\in\mathfrak{h}$, as the vector field $\zeta_Q \in \mathfrak{X}(Q)$ defined by
\begin{equation}
\zeta_Q(q) := \left.\frac{d}{d\varepsilon}\right|_{\varepsilon = 0}\Phi_{\exp_{H}( \varepsilon\zeta)}(q). \nonumber
\end{equation}
 The tangent and cotangent lifts $T\Phi$ and $T^*\Phi$ of $\Phi$ are group actions on tangent and cotangent spaces,
\begin{align}
&T\Phi: H \times TQ \rightarrow TQ, \quad (h, v_q) \mapsto T\Phi_h(v_q), \nonumber
\\ &T^*\Phi: H \times T^*Q \rightarrow T^*Q, \quad (h, \alpha_q) \mapsto (T\Phi_{h^{-1}})^*(\alpha_q). 
\label{cotan_lift_action}
\end{align}
The cotangent lift momentum map $\mathbf{J} : T^*Q \rightarrow \mathfrak{h}^*$ associated with the action $\Phi$ is determined, for arbitrary $\zeta \in \mathfrak{h}$, by
\begin{equation}\label{Definition_momentum_map}
\left<\mathbf{J} (\cdot), \zeta\right>_{\mathfrak{h}^* \times \mathfrak{h}}=\left<\,\cdot\,, \zeta_Q(q)\right>_{T^*Q \times TQ}= \theta_Q(\zeta_{T^*Q} (\cdot)),
\end{equation}
where $\theta_Q$
 is the canonical one form on $T^*Q$ and the vector field $\zeta_{T^*Q}\in \mathfrak{X}(T^*Q)$ is the infinitesimal generator for the cotangent lifted action \eqref{cotan_lift_action}. These cotangent lift momentum maps are an important special case of the general concept of a momentum map of a group acting on a Poisson manifold \cite{MaRa03}.

\paragraph{Noether's theorem.} 
We proceed with a version of Noether's theorem for systems with symmetry. The action $\Phi$ of a Lie group $H$ on $G$  induces an action on the space of curves $\mathcal{C}(G)$ and thus on both $TG$ and $T^{(2)}G$. Let us consider a Lagrangian $L: T^{(2)}G \to \mathbb{R}$, with the symmetry
\begin{equation}
L\left([\Phi_h\circ g]^{(2)}_{\Phi_h( g_0)}\right)=L\left([g]_{g_0}^{(2)}\right), \label{Lagrangian_invariance}
\end{equation}
 for all $ h\in H$ and any $\left[g\right]_{g_0}^{(2)} \in T^{(2)}G$. For consistency we continue to assume the Lagrangian is left-invariant, however the assumption is not essential for this section. For a curve $g(t) \in G$ consider the variation induced by the action of $H$,  $\delta_\zeta g = \zeta_G (g)$ for a $\zeta \in \mathfrak{h}$. Then
\begin{equation}\label{Symmetry_zeta}
\delta_\zeta \int_0^a L([g]) = 0.
\end{equation}
If $g(t)$ is a solution to Hamilton's principle \eqref{Hamiltons_principle}, then by Theorem \ref{Theorem_HOHP} there exists a lift $\gamma(t) \in P^{(2)}G$ satisfying the HOHP equations \eqref{RHP_dynamics}-\eqref{RHP_reconstruction}. By  \eqref{One-form_pull_back} and \eqref{Symmetry_zeta},
\begin{align} 
0 =  ds( \delta_\zeta x) = \theta(\delta_\zeta x(a)) - \theta(\delta_\zeta x(0))
= \theta(\zeta_{T^*(TG)}(x(a))) - \theta(\zeta_{T^*(TG)}( x(0))) , \nonumber
\end{align}
where $x(t) = \phi(\gamma(t)) \in T^*(TG)$. Using \eqref{Definition_momentum_map} we identify the cotangent lift momentum map $\mathbf{J}: T^*(TG) \to \mathfrak{h}^*$ associated with the action $T\Phi: H \times TG \to TG$ and find the conservation law
\begin{equation}
\mathbf{J}(x(t))=\mathbf{J}(x(0)) \qquad \text{for all} \ t\in[0,a].\nonumber
\end{equation}
\paragraph{Example.} Let $H = G$ and consider the left action $L$ of $G$ on itself,
\begin{equation}
  L: G \times G \rightarrow G, \quad L_h(g) = hg. \nonumber
\end{equation}
This is by definition a symmetry for any left-invariant Lagrangian.
The tangent lift $TL$ is given by
\begin{align}
  TL: G \times TG, \quad TL_h(g, \xi) =  (hg, \xi),\nonumber
\end{align}
where as usual we used left trivialisation, $TG \cong G \times \mathfrak{g}$. 
The cotangent lift momentum map associated with $TL$ is
\begin{equation} 
  \mathbf{J}: T^*(TG) \rightarrow \mathfrak{g}^*, \quad   (g, \xi, \mu, \nu) \mapsto  \operatorname{Ad}_{g^{-1}}^* \mu. \nonumber
\end{equation}
Therefore, $\mathbf{J}$ is conserved along solutions of the HOHP equations \eqref{RHP_dynamics}-\eqref{RHP_reconstruction}.

The Hamiltonian approach to Noether's theorem is discussed in \cite{MaRa03}. It is however the variational viewpoint that will expedite similar results for the variational integrators we develop next.
\section{Discrete-time mechanics}\label{discretetime}
Continuing with the Hamilton-Pontryagin approach of the previous section, we now address the issue of geometric discretisation. We proceed systematically with a discrete variational principle and, by considering end point terms, we establish the symplectic and momentum-preserving properties of the resulting numerical scheme. 
\subsection{Runge--Kutta--Munthe-Kaas (RKMK) methods}\label{Subsec-RKMK}
In this section we obtain a discretisation of the reconstruction equations \eqref{RHP_reconstruction}. This will be achieved by combining a standard Runge--Kutta (RK) method to integrate $\dot{\xi} = u$ with a Runge--Kutta--Munthe-Kaas (RKMK) method to integrate $\dot{g} = g\xi$. We shall introduce both types of methods in sufficient detail for our purposes here and refer to \cite{HLW06, IMKNZ00} for more complete discussions. 
\paragraph{Runge--Kutta (RK) method.} An $s$-stage RK method with step size $h$ may be employed to numerically integrate first-order ordinary differential equations of the form $\dot{y} = f(y, t)$ with $y(t_0) = y_0$. The task is to find solutions $Y^i$ for $i = 1, \ldots s$ that satisfy
\begin{align} \nonumber
  Y^i = y_0 +  h \sum_{j = 1}^s a_{ij}f(Y^j, t_0 + c_j h)
  \,.
\end{align}
The numerical estimate of $y(t_0 + h)$ is then given by 
\begin{align} \nonumber
  y_1 = y_0 + h \sum_{i = 1}^s b_i f(Y^i, t_0 + c_ih).
\end{align}
The $Y^i$ are known as \emph{internal stage variables} and
the real coefficients can be presented in a \emph{Butcher tableau}, 
\[
\begin{array}{c|ccc}c_1 & a_{11} & ... & a_{1s} \\ \vdots & \vdots &  & \vdots \\ c_s & a_{s1} & \ldots & a_{ss} \\\hline  & b_1 & \ldots & b_s\end{array}
\]
These coefficients must obey $c_i=\sum_ja_{ij}$ and $\sum_i b_i=1$ for an RK method of at least first order accuracy. A discussion of order conditions can be found in \cite{HLW06}.

 An $s$-stage RK method for the reconstruction equation $\dot{\xi} = f(\xi, t) = u(t)$ starting from $\xi(t_0) = \xi_0$ is thus
\begin{align}
  \Xi^i = \xi_0 + h\sum_{j = 1}^s a_{ij}V^j, \quad 
  \xi_1 = \xi_0 + h\sum_{i = 1}^s b_i V^i, \label{RKu}
\end{align}
where $V^j := u(t_0 + c_j h)$ and $\xi_1$ is our estimate of $\xi(t_0 + h)$.

\paragraph{Runge--Kutta--Munthe-Kaas (RKMK) method.} We now turn to the reconstruction equation $\dot{g} = g\xi$. In the RKMK method, one begins by rewriting the reconstruction equation on $\mathfrak{g}$, to which one then applies an RK method. Consider an approximation to the Lie exponential $\tau: \mathfrak{g} \rightarrow G$. Let $\tau$ be an analytic diffeomorphism in a neighbourhood of $0$ with $\tau(0) = e$ and $\tau(-\eta)\tau(\eta) = e$ for all $\eta \in \mathfrak{g}$. We define the right-trivialised differential $d\tau_\xi$ at $\xi \in \mathfrak{g}$ as,
\begin{equation}
  d\tau_\xi: \mathfrak{g} \rightarrow \mathfrak{g}, \quad \eta \mapsto (T_\xi\tau(\eta)) \tau(\xi)^{-1}, \nonumber
\end{equation}
with inverse $d\tau_\xi^{-1}$. We record the formula
\begin{equation}\label{Ad_tau_formula}
  d\tau_\xi^{-1} (\eta) = d\tau^{-1}_{-\xi} (\operatorname{Ad}_{\tau(-\xi)}\eta),
\end{equation}
which can be obtained by differentiating $\tau(-\xi)\tau(\xi) = e$. For small times $t$ a solution of $\dot{g} = g\xi$ with $g(0) = g_0$ can be parameterized as $g(t) = g_0 \tau(\Theta(t))$ for a $\Theta\in\mathcal{C}(\mathfrak{g})$. Then
\begin{equation} \nonumber
  g\xi = g_0 \tau(\Theta) \xi = \dot{g} = g_0 (d\tau_{\Theta}(\dot{\Theta})) \tau(\Theta),
\end{equation}
and therefore, using \eqref{Ad_tau_formula}, 
\[
\dot{\Theta} = f(\Theta, t) = d\tau^{-1}_{-\Theta} \xi(t)
\,.
\]
An RK method for this equation, with starting point $\Theta(t_0) = 0$, is
\begin{align}
  \Theta^i = \sum_{j = 1}^s a_{ij}d\tau^{-1}_{-h\Theta^j} \Xi^{j}, \quad
\Theta_1 = h\sum_{i = 1}^s b_i d\tau^{-1}_{-h\Theta^i}\Xi^i, \label{RKtheta}
\end{align}
where $\Xi^j:=\xi(t_0+c_jh)$. Define also $G^i = g_0 \tau(h\Theta^i)$. Equations \eqref{RKtheta} together with
\begin{equation}\label{RKg_step}
  g_1 = g_0 \tau(\Theta_1)
\end{equation}
constitute an $s$-stage RKMK method with step size $h$. Equation \eqref{RKg_step} is our estimate for $g(t_0 + h)$.  
The following theorem can be found in \cite{HLW06}.
{\em
\begin{theorem} \label{Thm_HLW}
Assume the underlying RK method in RKMK has order $p$ and $\tau$ is a $p^{th}$-order approximation to the Lie exponential. If one truncates $d\tau^{-1}$ at order $p-2$, the RKMK method is still of order $p$.
\end{theorem}
} 
For $p=2$, therefore, equations \eqref{RKtheta} can be simplified by approximating $d\tau^{-1}$ by the \emph{identity map}. The resulting RKMK method is of second order,
\begin{align}
\tau^{-1}(g_0^{-1}G^i) = h\sum_{j = 1}^s a_{ij}\Xi^{j}, \ \ \tau(g_0^{-1}g_1) = h\sum_{i = 1}^s b_i\Xi^i. \label{RKtheta_step_simplified}
\end{align}

\paragraph{Discrete reconstruction equations.} Our final method for effectively solving $(g^{-1}\dot{g})\dot{}=u$ comes in the two stage form achieved by combining the above methods, vector space RK and RKMK. In \eqref{RKtheta_step_simplified}, $\Xi^i$ are \emph{exactly} $\xi(t_0+c_i)$; in \eqref{RKu} the  $\Xi^i$ are thought of as being \emph{approximately} $\xi(t_0+c_i)$. The combined method feeds the calculated $\Xi^i$ in \eqref{RKu} into \eqref{RKtheta_step_simplified}. We now proceed by inserting the discrete reconstruction equations into a discrete HOHP principle.

\subsection{Higher-order Hamilton--Pontryagin Integrators}\label{Sec-HOHPI}
\paragraph{Discrete path space.} Our goal is to approximate the reduced higher-order Hamilton--Pontryagin principle (HOHP) in \eqref{HOHPP}. This is achieved by replacing the integral with a sum over discrete time points $t_k = kh$ for $k=0,\ldots, N$, where $h$ is the time step. The continuous variables are replaced by the $s$-stage RKMK variables discussed above and the kinematic constraints by the RKMK equations. A tentative choice of discrete path space is therefore given by the set of maps
 \begin{equation}
   \left\{\gamma: \{t_k\}_{k=0}^N \rightarrow T^*(TG) \times \left(P^{(2)}G\right)^s\right\} \nonumber
 \end{equation}
 which, by using the variables
\begin{align}
\gamma(t_k) =  \Big(g_k, \xi_k, \mu_k, \nu_k, \big\{G_i^k, \Xi_k^i, V_k^i, \mu_k^i, \nu_k^i\big\}_{s = 1}^N\Big), \nonumber
\end{align}
can be identified with the set $T^*(TG)^{N+1} \times \left(P^{(2)}G\right)^{s \cdot (N+1)}$.
 The discrete cost functional that we will define below in \eqref{Discrete_HOHP_cost_functional} will not depend explicitly on $\mu_0, \nu_0$ or $G_N^i, \Xi_N^i, V_N^i, \mu_N^i, \nu_N^i$. We therefore choose to omit these elements and define the \emph{discrete path space} $\mathcal{C}_d$ as follows,
\begin{align}
  \mathcal{C}_d &:= TG \times (T^*(TG))^N \times (P^{(2)}G)^{N\cdot s} \nonumber\\
&\cong (G\times \mathfrak{g}) \times (G \times \mathfrak{g} \times 2\mathfrak{g}^*)^N \times (G \times 2\mathfrak{g} \times 2\mathfrak{g}^*)^{N\cdot s}. \nonumber
\end{align}
We shall use the following index notation to denote elements,
\begin{equation}
  \left\{(g_0, \xi_0), \left(g_k, \xi_k, \check{\mu}_k, \nu_k\right)_{k=1}^N, \left(G_k^i, \Xi_k^i, V_k^i, \mu_k^i, \nu_k^i\right)_{i, k = 1, 0}^{s, N-1}\right\}. 
  \nonumber
\end{equation}
\paragraph{Discrete variational principle.} We now give the discrete HOHP integral. Define the functional $S: \mathcal{C}_d \rightarrow \mathbb{R}$ as
\begin{align}
\begin{split}
  S = h \sum_{k = 0}^{N-1} &\Bigg[\sum_{i=1}^s b_i \ell(\Xi_k^i, V_k^i) 
  + \Big<\mu_k^i, \frac{1}{h}\tau^{-1}(g_k^{-1}G_k^i) - \sum_{j=1}^s a_{ij}\Xi^j_k\Big>\\ & + \Big<\nu^i_k, \frac{1}{h}(\Xi_k^i - \xi_k) - \sum_{j=1}^s a_{ij}V_k^j\Big>\Bigg]+ \Big<\check{\mu}_{k+1}, \frac{1}{h} \tau^{-1}(g_k^{-1}g_{k+1}) - \sum_{i = 1}^s b_i \Xi_k^i\Big> \\ &+ \Big<\nu_{k+1}, \frac{1}{h}(\xi_{k+1} - \xi_k) - \sum_{i = 1}^s b_iV^i_k\Big>.\label{Discrete_HOHP_cost_functional}
\end{split}
\end{align}
Let $\gamma \in \mathcal{C}_d$ and $\delta \gamma \in T_\gamma \mathcal{C}_d$. Also denote $\eta_k = g_k^{-1}\delta g_k$ and $\eta_k^i = (G^i_k)^{-1}\delta G_k^i$. Moreover, set $\Xi_{k+1} = \sum_{j = 1}^s b_j\Xi^k_j$ and $\Theta^i_k = \sum_{j=1}^s a_{ij}\Xi^j_k$. Furthermore, it will be useful to define
\begin{align}
  \nu_0 = \nu_1 + \sum_{i = 1}^s \nu_0^i, \qquad 
  \mu_0 = (d\tau^{-1}_{h\Xi_1})^* \check{\mu}_1 + \sum_{i = 1}^s \left[(d\tau_{h\Theta_0^i})^*\mu_0^i\right], \label{mu0-nu0}
\end{align}
as well as
\begin{equation}\label{Coordinate_transform}
  \mu_k = (d\tau^{-1}_{-h\Xi_k})^*\check{\mu}_k, \qquad k = 1, \ldots, N. \nonumber
\end{equation}
In the following equations, the index $k$ takes range $0, \ldots N-1$ unless explicitly stated otherwise.  Discrete variations of $S$ now give,
\begin{align}
  \delta S  &= h\,\sum_k \Big(\Big[\sum_i \Big<\frac{1}{h}(d\tau^{-1}_{-h\Theta_k^i})^* \mu^i_k, \eta^i_k\Big>  + \Big<b_i \frac{\delta \ell}{\delta \xi} + \frac{1}{h}\nu^i_k - \sum_j a_{ij}\mu^j_k - b_i \check{\mu}_{k+1}, \delta \Xi_k^i\Big> \nonumber \\
&+ \Big<b_i \frac{\delta \ell}{\delta u} - \sum_j a_{ji}\nu^j_k - b_i \nu_{k+1}, \delta V_k^i\Big> + \left<\delta \mu^i_k, \_ \right> + \left<\delta \nu_k^i, \_\right>\Big] + \left<\delta \check{\mu}_{k+1}, \_ \right> + \left<\delta \nu_{k+1}, \_\right>\Big) \nonumber \\
&+\
\sum_{k=1}^{N-1} \Big( \Big<\nu_k - \nu_{k+1} - \sum_i\nu^i_k, \delta \xi_k\Big>  + \Big<\mu_k -(d\tau^{-1}_{h\Xi_{k+1}})^* \check{\mu}_{k+1} -\sum_i (d\tau^{-1}_{h\Theta_k^i})^*\mu_k^i, \eta_k\Big> \Big)\nonumber \\
&+\left<\mu_N, \eta_N\right> + \left<\nu_N, \delta \xi_N\right> - \left<\mu_0, \eta_0\right> - \left<\nu_0, \delta \xi_0\right>. \label{Discrete_variations}
\end{align}
A curve $\gamma\in \mathcal{C}_d$ is a solution of the HOHP principle, if $\delta S = dS(\delta \gamma) = 0$ for all vectors $\delta \gamma$ with fixed boundary points $g_0, g_N$ and velocities $\xi_0, \xi_N$, that is $\eta_0 = \eta_N = \delta \xi_0 = \delta \xi_N = 0$. We set the notation $\frac{\delta \ell^i_k}{\delta \xi} = \frac{\delta \ell}{\delta \xi}(\Xi_k^i, V_k^i)$ and obtain the RKMK reconstruction equations
\begin{align}
\begin{split}
&\Xi_k^i=\xi_k+h\sum_{j=1}^sa_{ij}V_k^j\,, \qquad  G_k^i=g_k\tau\Big(h\sum_{j=1}^sa_{ij}\Xi_k^j\Big),\\
&\xi_{k+1}=\xi_k+h\sum_{j=1}^sb_jV_k^j\,,  \qquad  g_{k+1}=g_k\tau\Big(h\sum_{j=1}^sb_{j}\Xi_k^j\Big)\,,
\end{split} \label{discrete_reconstruction}
\end{align}
the auxiliary equations
 \begin{align}
   \mu_k^i = 0\,,\qquad \nu_{k+1} = \nu_k - \sum_i\nu^i_k, \label{auxiliary}
 \end{align}
the discrete equations for the Ostrogradsky momenta,
\begin{align}
\frac{\delta \ell^i_k}{\delta \xi}=\left(d\tau_{-h\Xi_{k+1}}\right)^*\mu_{k+1} - \frac{1}{hb_i}\nu_k^i \,, \qquad
\frac{\delta \ell^i_k}{\delta u} = \sum_j \left(\frac{a_{ji}}{b_i}-1\right)\nu_k^j - \nu_k, \label{Ostrogradsky}
\end{align}
and the discrete version of the Euler--Poincar\'e equation,
\begin{equation}\label{Discrete_EP}
  \mu_{k+1} = (d\tau^{-1}_{-h\Xi_{k+1}})^* (d\tau_{h\Xi_{k+1}})^* \mu_k.
\end{equation}
\subsubsection{Geometric properties}
\paragraph{Symplecticity.} A solution $\gamma \in \mathcal{C}_d$ is said to have \emph{initial conditions} $(g_0, \xi_0, \mu_0, \nu_0) \in T^*(TG)$, upon using the definitions of $\nu_0$ and $\mu_0$ given in \eqref{mu0-nu0}. Assume the Lagrangian $\ell$ allows the internal variables $G_k^i,\Xi_k^i,V_k^i,\mu_k^i,\nu_k^i$ to be eliminated from the equations, which subsequently induce a discrete flow map $F: T^*(TG)\rightarrow T^*(TG)$. That is, for any solution $\gamma$ the flow map satisfies $F^k(g_0, \xi_0, \mu_0, \nu_0) = (g_k, \xi_k, \mu_k, \nu_k)$. We express the discrete HOHP cost functional in terms of initial conditions, as follows: Let $s: T^*(TG) \rightarrow \mathbb{R}$ be defined by $s (z) = S(\gamma)$, where $\gamma$ is the solution curve with initial conditions $z \in T^*(TG)$. A variation $\delta z \in T_zT^*(TG)$ corresponds to a variation $\delta \gamma \in T_\gamma\mathcal{C}_d$, and we compute the
generating-function relation,
\begin{align}
  ds(\delta z) =  \left<\mu_N, \eta_N\right> + \left<\nu_N, \delta \xi_N\right> - \left<\mu_0, \eta_0\right> - \left<\nu_0, \delta \xi_0\right> = \left[(F^*)^N\theta - \theta\right](\delta z) \label{discrete_end_points}
\,.\end{align}
 Upon taking the exterior derivative, one finds that the symplectic form $\omega$ is preserved by $F^N$. In particular, setting $N=1$ shows that the discrete flow map $F$ is symplectic.

\paragraph{Momentum preservation.} The argument continues along the same lines as in Section \ref{Sec-moma_noether}. Let us assume we have a (left-invariant) Lagrangian $L:T^{(2)}G\to\mathbb{R}$ and an action $\Phi$ of a Lie group $H$ on $G$. Note $(G_k^i,\Xi_k^i,V_k^i)\in G\times2\mathfrak{g}\cong T^{(2)}G$. For the group action to be a symmetry of the discrete dynamics we require that the Lagrangian is group invariant, as in \eqref{Lagrangian_invariance}, \emph{and} the group induced variations on $(G_k^i,\Xi_k^i,V_k^i)\in T^{(2)}G$ preserve the reconstruction relations \eqref{discrete_reconstruction}. That is, for $\zeta \in \mathfrak{h}$ and  $h(\varepsilon) := \exp_H(\varepsilon \zeta)$ the induced variation $\left(G_k^i(\varepsilon),\Xi_k^i(\varepsilon),V_k^i(\varepsilon)\right)$ still allow \eqref{discrete_reconstruction} to hold for some $\left\{g_k(\varepsilon),\xi_k (\varepsilon)\right\}_{k=0}^N$.  Similar to the continuous case, if all the underlying RK relations between the variables in $\mathcal{C}_d$ are implicitly assumed, then a solution  $(G_k^i,\Xi_k^i,V_k^i)$ to Hamilton's principle
\begin{equation}
\delta \sum_{k,i} b_i\ell(\Xi_k^i,V_k^i) h=0\,,
\end{equation}
can be lifted to a solution $\gamma$ of the discrete HOHP principle. Note that here variations satisfy fixed end point conditions $\delta g_0 = \delta g_N = 0$ and $\delta \xi_0 = \delta \xi_N = 0$. The group symmetry implies
\begin{equation}
\delta_{\zeta}\sum_{k,i} b_i \ell(\Xi_k^i,V_k^i) h=  0\,.
\end{equation}
For a variation $\delta_\zeta\gamma$ respecting \eqref{discrete_reconstruction} this implies $\delta_\zeta S = 0$, for $S$ as in \eqref{Discrete_HOHP_cost_functional}.
 From \eqref{discrete_end_points} we can then conclude that the cotangent-lift momentum map $\mathbf{J}: T^*(TG) \to \mathfrak{h}^*$ associated with the action $T\Phi: H \times TG \to TG$ is preserved by the HOHP integration scheme,
\begin{align} 
0 =  ds( \delta_\zeta z)
= \theta(\zeta_{T^*(TG)}(z_N)) - \theta(\zeta_{T^*(TG)}( z_0)) = \mathbf{J}(z_N)-\mathbf{J}(z_0), \nonumber
\end{align}
where we wrote $z_0 \in T^*(TG)$ for the initial condition of $\gamma$ and $z_N := F^N(z_0)$.
 \paragraph{Example.}  As in the continous-time case, let $H = G$ and consider the cotangent-lift momentum map associated with $TL: G \times TG \to TG$, the tangent lift of the left action of $G$ on itself. We recall that
\begin{equation}
  \mathbf{J}: T^*(TG) \rightarrow \mathfrak{g}^*, \quad   (g, \xi, \mu, \nu) \mapsto  \operatorname{Ad}_{g^{-1}}^* \mu. \nonumber
\end{equation}
By the above discussion, $\mathbf{J}$ is conserved along solutions of the discrete HOHP equations \eqref{discrete_reconstruction}--\eqref{Discrete_EP}. That is, $\mathbf{J}_k := \operatorname{Ad}^*_{g_k^{-1}}\mu_k$ is conserved, $\mathbf{J}_{k+1} = \mathbf{J}_k$. 

We remark that this conservation law can alternatively be obtained directly from the HOHP equations and a property therein which relates them to the Euler--Poincar\'e equations. It follows from \eqref{Ad_tau_formula} that for any $\eta, \xi \in \mathfrak{g}$ and $\mu \in \mathfrak{g}^*$, 
\begin{align*}
\left< \left(d\tau_{-\xi}^{-1}\right)^*(\mu), \eta\right> =\left<\mu,  d\tau_{-\xi}^{-1} (\eta)\right> = \left<\mu,  d\tau^{-1}_{\xi} (\operatorname{Ad}_{\tau(\xi)}\eta)\right> 
= \left<\operatorname{Ad}^*_{\tau(\xi)}(d\tau^{-1}_{\xi})^*(\mu), \eta\right>.
\end{align*}
Therefore
\begin{align}\label{Ad_star_equation}
  \operatorname{Ad}^*_{\tau(\xi)} = (d\tau^{-1}_{-\xi})^* \circ (d\tau_\xi)^*. 
\end{align}
For a solution of the discrete HOHP equations, \eqref{Discrete_EP} implies
\begin{align*}
  \mu_{k+1}& = (d\tau_{-h\Xi_{k+1}}^{-1})^* (d\tau_{h\Xi_{k+1}})^* \mu_k = \operatorname{Ad}^*_{\tau(h\Xi_{k+1})}\mu_k \\ &= \operatorname{Ad}^*_{\tau(h\Xi_{k+1})} \operatorname{Ad}^*_{g_k} \operatorname{Ad}^*_{g_k^{-1}} \mu_k =  \operatorname{Ad}^*_{g_{k+1}} \operatorname{Ad}^*_{g_k^{-1}} \mu_k.
\end{align*}
Hence,
\begin{align*}
  \operatorname{Ad}_{g_{k+1}^{-1}}^* \mu_{k+1} =   \operatorname{Ad}_{g_{k}^{-1}}^* \mu_{k} ,
\end{align*}
which is the conservation law for the discrete momentum, $\mathbf{J}_{k+1} = \mathbf{J}_k$. 

\subsubsection{Implementation}\label{AppCubics}
 Our primary application of the HOHP scheme addresses the Riemannian cubics on $SO(3)$ with a bi-invariant metric. We implement initial value problem solvers using the implicit Euler method and the St\"{o}rmer--Verlet method.
 \begin{figure}
\centering
\includegraphics[scale=0.3]{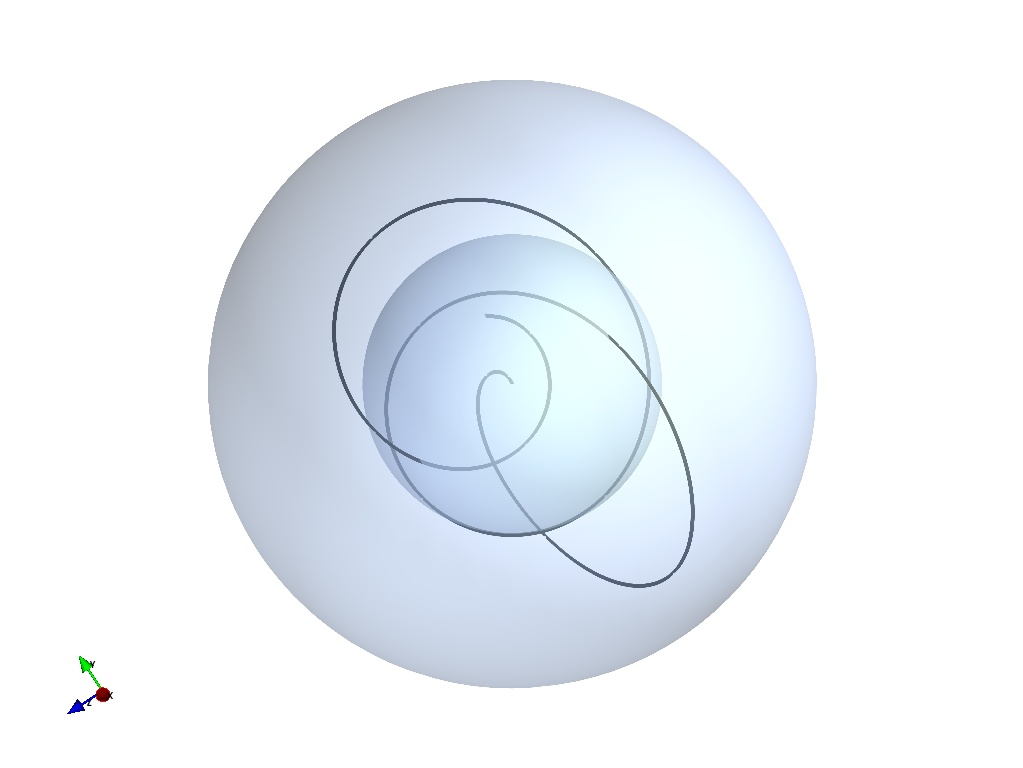}
\caption{\footnotesize An illustration of the  St\"{o}rmer--Verlet method applied to the initial value problem of Riemannian cubics on the rotation group, presented in \eqref{dep}--\eqref{SVxi}. In this figure, a given element of the group corresponds to a point on the radial line along the rotation axis. The distance of the point from the origin represents the rotation angle measured in a right-handed coordinate system. The center and the boundary of the outer sphere therefore both represent the identity matrix. The displayed curve results from numerically integrating system \eqref{dep}--\eqref{SVxi} over a period of $2\pi$. The chosen initial conditions were  $g_0=\Id$, $\xi_0=(-6,1,0)$, $\nu_0=(0,0,6)$ and $\mu_0=(0,36,0)$, which belong to certain $2\pi$-periodic solutions found in \cite{No2004}.}
 \label{SV_IVP_plot}
\end{figure}
\paragraph{Numerical schemes.} We consider two instances of the general RKMK scheme, an implicit Euler method and a St\"{o}rmer--Verlet type method given respectively by,
\[
\begin{tabular}{c|c}1 & 1 \\\hline  & 1\end{tabular} \ \ \qquad\text{and} \ \ \qquad \begin{tabular}{c|cc}0 & 0 & 0 \\ 1 & 1/2 & 1/2 \\\hline  & 1/2 & 1/2\end{tabular}
\]
For the Riemannian cubic Lagrangian, $\ell(\xi,u)=\frac{1}{2}\|u\|^2$, the Euler method gives
\begin{align}
\mu_{k+1}&=(d\tau^{-1}_{-h\xi_{k+1}})^*(d\tau_{h\xi_{k+1}})^*\mu_k\,,\nonumber \\
\nu_{k+1}&=\nu_k-h(d\tau_{-h\xi_{k+1}})^*\mu_{k+1}\,,\nonumber \\
g_{k+1}=g&_k\tau(h\xi_{k+1})\,, \quad \xi_{k+1}=\xi_k+h\nu_k^\sharp\,. \nonumber
\end{align}
This method is explicit and straight-forward to implement. For St\"{o}rmer--Verlet the result is similar but implicit 
\begin{align}
\mu_{k+1}=(d\tau^{-1}_{-h\Xi_{k+1}})^*&(d\tau_{h\Xi_{k+1}})^*\mu_k \,, \quad \Xi_{k+1}=\frac{1}{2}\left(\xi_{k+1}+\xi_k\right), \label{dep}\\
\nu_{k+1}&=\nu_k - h(d\tau_{h\Xi_{k+1}})^*\mu_k\,,\\
g_{k+1}=g_k\tau(h\Xi_{k+1})& \,, \quad \xi_{k+1}=\xi_k+h \left(\nu_k-\frac{h}{2}(d\tau_{h\Xi_{k+1}})^*\mu_k\right)^\sharp\,. \label{SVxi}
\end{align}
The most notable differences are the appearance of the average velocity $\Xi_{k+1}$ and the second order approximation in the $\xi_{k+1}$ equation. To implement this method we can apply a fixed point method to find $\xi_{k+1}$ with the function
\begin{align}
f_k(\xi)=\xi-\xi_k-h\left(\nu_k+\frac{h}{2}\left(d\tau_{\frac{h}{2}(\xi+\xi_k)}\right)^*\mu_k\right)^\sharp\,. \nonumber
\end{align}
 An example of a map $\tau$ is the Cayley map, $\cay:\mathfrak{so}(3)\to SO(3)$, given by
\begin{align}
\cay(\widehat{X})=\left(e-\widehat{X}/2\right)^{-1}\left(e+\widehat{X}/2\right), \quad
d\cay_{\widehat{X}}\widehat{Y} = \left(e-\widehat{X}/2\right)^{-1}\widehat{Y}\left(e+\widehat{X}/2\right)^{-1}
\end{align}
where $e$ is the identity matrix. With the Cayley map we can implement the St\"{o}rmer--Verlet method, say, on $SO(3)$ to produce plots such as in Figure \ref{SV_IVP_plot}. Note that both algorithms give discrete flow maps $F:T^*TSO(3)\to T^*TSO(3)$ which are symplectic and momentum preserving, as discussed in Section \ref{Sec-HOHPI}. 
\section{Application to  trajectory planning}\label{applTM-sec}
In this section we discuss an application of the discrete HOHP variational principle to a second order trajectory-planning problem, where the interpolating curve evolves by means of a group action. One aims at finding a curve that passes near a series of given target points at prescribed times. A precise definition of the problem will be given below. An earlier investigation of the same problem in \cite{Gay-BalmazEtAl2010} found that optimal curves satisfy second order Euler--Poincar\'e equations between time nodes. At the time nodes the (otherwise conserved) momentum experiences a kick that is related in a simple fashion with the optimal mismatch between interpolating curve and target at the respective node, leading to a piecewise constant momentum. The goal of the present section is to show that this geometric characteristic carries over to a discrete-time algorithm.
\subsection{Continuous-time trajectory planning}
\paragraph{Background.} Our motivation to study this type of trajectory-planning problem lies with potential applications in computational anatomy, the modeling and quantifying of diffeomorphic shape evolution \cite{MillerARBE2002,MillerIJCV2001}. Usually one seeks a geodesic path on the space of shapes between given initial and final data. This approach can be adapted for longitudinal data interpolation, interpolation through a sequence of data points, by piecewise geodesic curves. A class of higher-order generalizations providing smoother solution curves than the piecewise geodesic ones was studied in \cite{Gay-BalmazEtAl2010} in the finite dimensional setting.

Here we first reformulate the problem in terms of the HOHP principle, which is amenable to the discretisation introduced in Section \ref{Sec-HOHPI}.
\paragraph{Problem formulation.} Consider a left representation of the Lie group $G$ on a vector space $V$ with norm $\|.\|_V$,
\begin{equation}
  G \times V \to V, \qquad (g, I) \mapsto gI. \nonumber
\end{equation}
The problem under consideration is defined as follows: 
\textit{Given a Lagrangian $\ell: 2 \mathfrak{g} \rightarrow \mathbb{R}$, $\sigma, t_1, \ldots, t_l \in \mathbb{R}$, $T_0, I_{t_1}, \ldots , I_{t_l} \in V$, and $\xi^0_0 \in \mathfrak{g}$, minimise the functional
\begin{equation}\label{LongitudinalInterpolation}
S := \int_0^{t_l} \ell(\xi,u) + \left<\mu, g^{-1} \dot{g} - \xi\right> + \left<\nu, \dot{\xi} - u\right>  dt+ \frac{1}{2\sigma^2}\sum_{i = 1}^l \left\|g^{-1}(t_i)T_0 - I_{t_i}\right\|^2_V
\end{equation}
subject to the conditions $\xi(0) = \xi^0_0$ and $g(0) = e$. Variations are taken amongst curves $(g, \xi, u, \mu, \nu)(t) \in P^{(2)}G$, where we require that $g \in C^1([0, 1])$ and $g|_{[t_i, t_{i+1}]} \in C^\infty([t_i, t_{i+1}])$.}

This type of trajectory-planning problem is familiar, for example, from computational anatomy, where one would typically think of $T_0$ as a template that is deformed by a curve of diffeomorphisms $g^{-1}(t)$, in turn generated by the time-dependent vector field $\xi(t)$. At times $t_i$ the curve passes near the given targets $I_{t_i}$, the parameter $\sigma$ determining the proximity of the passage. In that case the Lie group $G$ is infinite dimensional. It is clear that this type of trajectory-planning problem is also relevant in numerous finite-dimensional situations  whenever a path $g(t) T_0$ generated by transformations $g(t)$ is required to pass by prescribed target points $I_{t_i}$ at given times $t_i$.  Here we focus on the latter case.

\paragraph{Euler-Lagrange equations.}  We assume that the norm on $V$ is induced by an inner product $\left\langle\,\cdot\,,\,\cdot\,\right\rangle_V$ and define the isomorphism $\flat: V \to V^*$ by $\left<u, v\right>_V = \left<u^\flat, v\right>$. We denote the cotangent-lift momentum map associated with the action of $G$ on $V$ by $\diamond: V \times V^* \rightarrow \mathfrak{g}^*$, that is (see \eqref{Definition_momentum_map}),
\begin{equation}
  \left\langle I\diamond \omega, \xi \right\rangle_{\mathfrak{g}^* \times \mathfrak{g}} = \left\langle \omega, \xi_V (I)\right\rangle_{V^* \times V}, \quad \text{for all} \quad I \in V, \; \omega \in V^*,  \;\xi \in \mathfrak{g}.
\label{diamond_defn}
\end{equation}
We compute the Euler-Lagrange equations governing any minimiser of $S$ in analogy with earlier calculations done in Section \ref{Sec-Variational_Principles}.
The only complication arises from taking variations of the penalty term in \eqref{LongitudinalInterpolation}. Setting $\eta := g^{-1} \delta g$, we obtain
\begin{align*}
  \delta \frac{1}{2\sigma^2}\sum_{i = 1}^l \left\|g(t_i)^{-1}T_0 - I_{t_i}\right\|^2_V &= -\frac{1}{\sigma^2} \sum_{i=1}^l \left<\left(g(t_i)^{-1}T_0 - I_{t_i}\right)^\flat, \left(\eta(t_i)_V \right)\left(g(t_i)T_0\right)\right>_{V^*\times V} \\
&=- \frac{1}{\sigma^2} \sum_{i=1}^l \left< \left(g(t_i)^{-1}T_0\right) \diamond \left(g(t_i)^{-1}T_0 - I_{t_i}\right)^\flat, \eta(t_i)\right>_{\mathfrak{g} \times \mathfrak{g}^*}.
\end{align*}
Taking variations of \eqref{LongitudinalInterpolation} and requiring $\delta S = 0$ yields the HOHP equations \eqref{RHP_dynamics}--\eqref{RHP_reconstruction} on  open intervals $(t_i, t_{i+1})$. We recall these equations as
\begin{align}
\dot\mu=\ad^*_\xi\mu, \qquad \mu=\frac{\delta \ell}{\delta \xi}-\dot{\nu}, \qquad  \nu= \frac{\delta \ell}{\delta u}, \qquad 
\dot{g}=g\xi, \qquad  \dot{\xi}=u. \label{HOHP_recall}
\end{align}
Moreover, writing $t_i^\pm$ for the limits from above and below respectively,
\begin{align}
  &\nu(t_i^+) - \nu(t_i^-) = 0, \label{Continuous_nu}\\
& \mu(t_i^+) - \mu(t_i^-) + \frac{1}{\sigma^2}\left(g(t_i)^{-1}T_0\right) \diamond \left(g(t_i)^{-1}T_0 - I_{t_i}\right)^\flat = 0, \label{Kick_mu}
\end{align}
 for $i = 1, \ldots l-1$, and
\begin{align}
  &\nu(t_l) = 0, \label{nu_final} \\
&\mu(t_l) -  \frac{1}{\sigma^2}\left(g(t_l)^{-1}T_0\right) \diamond \left(g(t_l)^{-1}T_0 - I_{t_l}\right)^\flat = 0. \label{mu_final}
\end{align}
\paragraph{Momentum kicks.} Optimal curves satisfy the HOHP equations \eqref{HOHP_recall}  between time nodes. While the Ostrogradsky momentum $\nu$ is continuous in time \eqref{Continuous_nu}, the Ostrogradsky momentum $\mu$  experiences a discontinuity (kick) at the time nodes. The kick is related to the optimal mismatch between interpolating curve and target points at the respective node, \eqref{Kick_mu}. The final values at $t_l$ of $\nu$ and $\mu$ are given in \eqref{nu_final} and \eqref{mu_final}.
The prescription of target points $I_{t_i}$ breaks the symmetry of the problem, which is reflected in the fact that the momentum $\mathbf{J} := \operatorname{Ad^*_{g^{-1}} \mu}$ is no longer preserved. Indeed, $\mathbf{J}$ is piecewise constant, with discontinuities at the time nodes given by
\begin{align}\label{Continuous_momentum_kick}
  \mathbf{J}(t_i^+) - \mathbf{J}(t_i^-) =- \frac{1}{\sigma^2} \operatorname{Ad}^*_{g(t_i)^{-1}} \left(\left(g(t_i)^{-1} T_0\right) \diamond \left(g(t_i)^{-1} T_0 - I_{t_i}\right)^\flat \right).
\end{align}
\subsection{Discrete-time trajectory planning}
As in Section \ref{Sec-HOHPI} we  replace the integral in \eqref{LongitudinalInterpolation} with a sum over time points $t_k = kh$ for $k=0,\ldots, N$, where $h$ is the time step. Again the continuous variables are replaced by the $s$-stage RKMK variables and the kinematic constraints by the RKMK equations, see Section \ref{Subsec-RKMK}. Let the target times be $t_i = N_{i}h$ for $i = 1, \ldots l$ with $N_l = N$. For convenience we also define $N_0 := 0$. The discretisation of \eqref{LongitudinalInterpolation} is given by
\begin{align}
\begin{split}
  S = h \sum_{k = 0}^{N-1} &\Bigg[\sum_{i=1}^s b_i \ell(\Xi_k^i, V_k^i) 
  + \Big<\mu_k^i, \frac{1}{h}\tau^{-1}(g_k^{-1}G_k^i) - \sum_{j=1}^s a_{ij}\Xi^j_k\Big>\\ & + \Big<\nu^i_k, \frac{1}{h}(\Xi_k^i - \xi_k) - \sum_{j=1}^s a_{ij}V_k^j\Big>\Bigg]+ \Big<\check{\mu}_{k+1}, \frac{1}{h} \tau^{-1}(g_k^{-1}g_{k+1}) - \sum_{i = 1}^s b_i \Xi_k^i\Big> \\ &+ \Big<\nu_{k+1}, \frac{1}{h}(\xi_{k+1} - \xi_k) - \sum_{i = 1}^s b_iV^i_k\Big> +  \frac{1}{2\sigma^2}\sum_{i = 1}^l \left\|g_{N_i}^{-1}T_0 - I_{t_i}\right\|^2_V \label{Discrete_TM_functional}
\end{split}
\end{align}
\paragraph{Discrete Euler-Lagrange equations.}
For the variations of the penalty term we obtain
\begin{align}
  \delta \frac{1}{2\sigma^2}\sum_{i = 1}^l \left\|g_{N_i}^{-1}T_0 - I_{t_i}\right\|^2_V =- \frac{1}{\sigma^2} \sum_{i=1}^l \left< \left(g_{N_i}^{-1}T_0\right) \diamond \left(g_{N_i}^{-1}T_0 - I_{t_i}\right)^\flat, \eta_{N_i}\right>_{\mathfrak{g} \times \mathfrak{g}^*}, \nonumber
\end{align}
where we set $\eta_k := g_k^{-1}\delta g_k$. As in \eqref{Coordinate_transform} we define 
\begin{equation}
  \mu_k = (d\tau^{-1}_{-h\Xi_k})^*\check{\mu}_k, \qquad k = 1, \ldots, N. \label{coord_change}
\end{equation}
 Taking variations of \eqref{Discrete_TM_functional} and requiring $\delta S = 0$ gives, for $ k = 0, \ldots N-1$, 
\begin{align}
\begin{split}
&\Xi_k^i=\xi_k+h\sum_{j=1}^sa_{ij}V_k^j\,, \qquad  G_k^i=g_k\tau\Big(h\sum_{j=1}^sa_{ij}\Xi_k^j\Big), \\
&\xi_{k+1}=\xi_k+h\sum_{j=1}^sb_jV_k^j\,,  \qquad  g_{k+1}=g_k\tau\Big(h\sum_{j=1}^sb_{j}\Xi_k^j\Big)\,, \\
   &\mu_k^i = 0\,, \qquad \frac{\delta \ell^i_k}{\delta \xi}=\left(d\tau_{-h\Xi_{k+1}}\right)^*\mu_{k+1} - \frac{1}{hb_i}\nu_k^i \,, \qquad
\frac{\delta \ell^i_k}{\delta u} = \sum_j \left(\frac{a_{ji}}{b_i}-1\right)\nu_k^j - \nu_k, \\
 &\nu_{k+1} - \nu_k + \sum_i\nu^i_k = 0. \label{Discrete_HOHP_TM_equations}
\end{split}
 \end{align}
Moreover, for interior indices $k \neq N_i \, \, (i =0, \ldots l)$, 
\begin{equation}\label{mu_discrete_interior}
\mu_k -(d\tau^{-1}_{h\Xi_{k+1}})^* (d\tau_{-h\Xi_{k+1}})^* \mu_{k+1} = 0, 
\end{equation}
and for node indices $N_i \,\, (i = 1, \ldots l-1)$, 
\begin{equation}\label{mu_discrete_node}
\mu_{N_i} -(d\tau^{-1}_{h\Xi_{N_i+1}})^* (d\tau_{-h\Xi_{N_i+1}})^* \mu_{N_i+1} - \frac{1}{\sigma^2}   \left(g_{N_i}^{-1}T_0\right) \diamond \left(g_{N_i}^{-1}T_0 - I_{t_i}\right)^\flat = 0. 
\end{equation}
Finally,
\begin{align}
  &\nu_N = 0, \label{nu_discrete_endpoint} \\
  &\mu_N - \frac{1}{\sigma^2} \left(g_N^{-1}T_0\right) \diamond \left( g_N^{-1}T_0 - I_{t_l}\right)^\flat = 0. \label{mu_discrete_endpoint}
\end{align}
\paragraph{Discrete momentum kicks.} Note that equations \eqref{Discrete_HOHP_TM_equations} are a subset of the HOHP equations of Section \ref{Sec-HOHPI}. The last equation in  \eqref{Discrete_HOHP_TM_equations} is the discrete version of \eqref{Continuous_nu}. The discrete Euler--Poincar\'e equation \eqref{mu_discrete_interior} completes the discrete HOHP equations for interior indices. This corresponds to the continuous-time equations \eqref{HOHP_recall}, which hold on open intervals. At time nodes $N_i$, the update formula for $\mu_{N_i+1}$ acquires an extra term, described by \eqref{mu_discrete_node}. This is analogous to the discontinuity of $\mu(t)$ seen in \eqref{Kick_mu}. The exact correspondence between equations \eqref{nu_discrete_endpoint} and \eqref{mu_discrete_endpoint} with \eqref{nu_final} and \eqref{mu_final} is obvious. The kicks of the discrete momentum $\mathbf{J}_k := \operatorname{Ad}^*_{g_k^{-1}}\mu_k$ are computed from \eqref{mu_discrete_interior} using \eqref{Ad_star_equation},
\begin{equation} \label{Discrete_momentum_kick}
  \mathbf{J}_{N_i+1} - \mathbf{J}_{N_i}  =- \frac{1}{\sigma^2} \operatorname{Ad}^*_{g_{N_i}^{-1}}\left( \left(g_{N_i}^{-1} T_0\right) \diamond \left(g_{N_i}^{-1} T_0 - I_{t_i}\right)^\flat \right),
\end{equation}
for $ i = 1, \ldots l-1$. This is the exact analogue of \eqref{Continuous_momentum_kick}. If $k$ is an interior index, then $\mathbf{J}_{k+1} = \mathbf{J}_k$. This discussion shows that the geometric behaviour of continuous-time momentum carries over to the discrete-time setting described above.

\subsection{Implementation via shooting method}
In this section we address the practical problem of minimizing the discrete functional  \eqref{Discrete_TM_functional}. Since any minimizing curve is known to satisfy the shooting equations \eqref{Discrete_HOHP_TM_equations}--\eqref{mu_discrete_node}, the search can be restricted to the space of solutions to these equations which results in a problem of reduced dimensionality (see also \cite{CoRiViRu2011}). The functional \eqref{Discrete_TM_functional} is accordingly written as a function of initial momenta $\mu_0$ and $\nu_0$,  $\mathcal{J}: 2\mathfrak{g}^* \to \mathbb{R}$. A gradient descent algorithm is then employed to carry out the minimization of $\mathcal{J}$. We discuss in some detail an efficient way of computing the gradient $ \nabla \mathcal{J}$ and present numerical simulations for the case of $SO(3)$ acting on $\mathbb{R}^3$. 
\paragraph{Shooting method.} Explicitly,  $\mathcal{J}: 2\mathfrak{g}^* \to \mathbb{R}$ is
\begin{equation}\label{implicit_func}
\mathcal{J}(\mu_0,\nu_0) = \int_0^{t_l}  \ell(\xi, u) dt+ \frac{1}{2\sigma^2}\sum_{i = 1}^l \left\|g^{-1}(t_i)T_0 - I_{t_i}\right\|^2_V, 
\end{equation}
where equations \eqref{HOHP_recall}-\eqref{Kick_mu} are implied, with initial conditions $g(0) = e$, $\xi(0) = \xi_0^0$ (recall that these are prescribed by the interpolation problem), $\mu(0) = \mu_0$ and $\nu(0) = \nu_0$. That is to compute  $\mathcal{J}(\mu_0,\nu_0)$, integrate \eqref{HOHP_recall}-\eqref{Kick_mu} with the given initial conditions, then evaluate the right hand side of \eqref{implicit_func}. 

In discrete time one proceeds in analogous fashion, defining  $\mathcal{J}_d: 2 \mathfrak{g}^* \to\mathbb{R}$ as
\begin{equation}
\mathcal{J}_d(\mu_0,\nu_0) =  h \sum_{k = 0}^{N-1}\sum_{i=1}^s b_i \ell(\Xi_k^i, V_k^i) 
   + \frac{1}{2\sigma^2}\sum_{i = 1}^l \left\|g^{-1}_{N_i}T_0 - I_{t_i}\right\|^2_V\,,
\label{disc_implicit_func}
\end{equation}
 the discrete shooting equations \eqref{Discrete_HOHP_TM_equations}--\eqref{mu_discrete_node} being implied. We then employ a gradient descent method to minimise $\mathcal{J}_d$. By construction, the corresponding solution curve obeys momentum map conservation for interior indices and has momentum kicks \eqref{Discrete_momentum_kick} at node indices. 
\subsubsection{Computation of the gradient} 
We now discuss the computation of the gradient $\nabla\mathcal{J}_d$ via adjoint equations, the concept of which we first expound in continuous time before presenting the discrete case. For this discussion we choose the bi-invariant cubic Lagrangian $\ell(\xi,u)=\frac{1}{2}\|u\|^2_{\mathfrak{g}}$.
\paragraph{Continous time.} For simplicity, let $l=1$. Introduce adjoint variables $P^0,P^1,V^0,V^1,V^2\in 2\mathfrak{g}^*\times3\mathfrak{g}=:A$ used to explicitly enforce the implied constraints in the definition of $\mathcal{J}$. Take $\mathcal{S}:\mathcal{C}(P^{(2)}G\times A)\to\mathbb{R}$ to be
\begin{align}
\begin{split}
\mathcal{S}(\alpha) = &\int_0^{t_1}\frac{1}{2}\|u\|_{\mathfrak{g}}^2 dt+ \frac{1}{2\sigma^2} \left\|g^{-1}(t_1)T_0 - I_{t_1}\right\|^2_V \\
+&\int_0^{t_1}\left< P^0, g^{-1}\dot g - \xi\right>+\left< P^1, \dot\xi - \nu^\flat\right> + \left< \mu + \dot\nu, V^0\right> \\
&\qquad+ \left< \dot\mu-\ad_\xi^*\mu, V^1\right> + \left<\nu-u^\flat, V^2\right>\, dt\,. \label{cts_adjoint_func}
\end{split}
\end{align}
Note that if $\widetilde\alpha$, the projection of the curve $\alpha \in \mathcal{C}(T^*(TG)\times A)$ onto $T^*(TG)$, satisfies the HOHP cubic equations \eqref{HOHPcubic} with initial condition $\widetilde{\alpha}(0)= (g_0, \xi_0, \mu_0,\nu_0)$ then $\mathcal{S}(\alpha) = \mathcal{J}(\mu_0,\nu_0)$. Taking variations while keeping $(g(0),\xi(0))$ fixed (the prescribed part of the initial conditions) we find that if such a curve $\alpha$ obeys terminal conditions,
\begin{align}
\begin{split}
&P^0(t_1)=\frac{1}{\sigma^2} \left(g(t_1)^{-1}T_0\right) \diamond \left( g(t_1)^{-1}T_0 - I_{t_1}\right)^\flat\,,\\ 
 &P^1(t_1)=0 \,, \quad V^0(t_1)= 0\,, \quad V^1(t_1)=0\,, \label{terminal_conditions}
\end{split}
\end{align} 
and adjoint equations
\begin{align}
\begin{split}
&\dot{P^0} = \ad_\xi^*P^0\,, \quad \dot{P^1}=-P^0 + \ad^*_{V^1}\mu\,,  \\
&\dot{V^0}= u -(P^1)^\sharp \,,  \quad  
\dot{V^1}=V^0-\ad_\xi V^1\,,  \label{Backward_shooting_continuous}
\end{split}
\end{align}
then the gradient $\nabla \mathcal{J}$ is given by
  \begin{align}
\nabla_{\mu_0} \mathcal{J} &= -V^1(0), \qquad \nabla_{\nu_0}\mathcal{J} = - V^0(0). 
 \label{func_grad}
\end{align}
In practice, one integrates \eqref{HOHPcubic} forward and \eqref{Backward_shooting_continuous} backward in time, with initialisation at $t = t_1$ given by \eqref{terminal_conditions}.
\paragraph{Discrete time.} For the purposes of implementation we are primarily interested in the discrete case and we now derive adjoint equations consistent with our discretisation of the forward shooting equations. We now allow for $l \geq 1$, but for simplicity restrict to the Euler type HOHP integration scheme for which the equations with momentum kicks \eqref{Discrete_HOHP_TM_equations}--\eqref{mu_discrete_node} can be written as, for $k=0,\ldots,N-1$,
\begin{align}
\begin{split}
g_{k+1}&=g_k\tau(h\Xi_{k+1})\,, \quad \Xi_{k+1}=\xi_k+h u_k\,, \quad \xi_{k+1}=\xi_k+h u_k\,, \quad u_k = \nu_k^\sharp\,, \\
\quad \nu_{k+1}&=\nu_k-h\check{\mu}_{k+1}\,, \quad
(d\tau^{-1}_{h\xi_{k+1}})^*\check\mu_{k+1}=(d\tau^{-1}_{-h\xi_{k}})^*\check\mu_k  +  \Phi_k(g_k)\,, 
\end{split} \label{euler_HOHP_TM}
\end{align} 
where $\Phi_k(g) := 0$ for $k\neq N_i$, $i=1,\ldots, l$ and 
\begin{align}
\Phi_k(g):=  - \frac{1}{\sigma^2}\left(g^{-1}T_0\right) \diamond \left(g^{-1}T_0 - I_{t_i}\right)^\flat , \quad \text{for}\quad k=N_i, \quad i=1,\ldots,l \,.
\end{align}
The discrete version of \eqref{cts_adjoint_func} is
\begin{align}
\begin{split}
\mathcal{S}_d(\alpha_d) =  \sum_{k=0}^{N-1}&\frac{h}{2}\|u_k\|_{\mathfrak{g}}+\left< P^0_{k+1}, \tau^{-1}(g^{-1}_k\dot g_{k+1}) - h\Xi_{k+1}\right> + \left< P^1_{k+1}, \xi_{k+1}-\xi_k - hu_k\right> \\
&+ \left< P^2_{k+1}, \Xi_{k+1}- \xi_k -hu_k \right> + \left<\nu_{k+1}-\nu_k + h\check\mu_{k+1}, V^0_{k+1}\right> \\
&+ \left<(d\tau^{-1}_{h\xi_{k+1}})^*\check\mu_{k+1}-(d\tau^{-1}_{-h\xi_k})^*\check\mu_k - \Phi_k(g_k),V_{k+1}^1\right> + \left< u_k - \nu_k, V^2_{k+1}\right> \\
&\hspace{-15pt}+\frac{1}{2\sigma^2}\sum_{i = 1}^l \left\|g_{N_i}^{-1}T_0 - I_{t_i}\right\|^2_V\,.
\end{split}
\end{align}
Notice that the discrete adjoint space $A_d$ has an extra variable $P^2\in\mathfrak{g}^*$. 
The resulting discrete adjoint equations consist of the terminal conditions
\begin{align}
  \begin{split}
    &\qquad(d\tau^{-1}_{-h\xi_{{N}}})^*P^0_{{N}} = - \Phi_N(g_N), \\
& P^1_{{N}} = - D^+_{{N}}V^1_{{N}}= 0\,, \quad
  V^0_{N} = 0\,, \quad
 V^1_{{N}} =0\,,
  \end{split}
\end{align}
and the discrete adjoint equations
\begin{align}
\begin{split}
P^0_{k} &= (d\tau_{-h\xi_{k}})^*\left((d\tau^{-1}_{h\xi_{k+1}})^*P^0_{k+1} -  \mathcal{A}_k(g_k,V^1_{k+1})\right)\,,   \\
P^1_{k} &= hP^0_{k+1} + P^1_{k+1} + D^-_kV_{k+1}^1-D^+_{k}V^1_{k}\,,\\
V^0_{k} &= V^0_{k+1}+h(P_{k+1}^1+ hP^0_{k+1} - u_k)\,, \\
 V^1_{k} &= d\tau_{h\xi_{k}}\left(d\tau^{-1}_{-h\xi_k}V^1_{k+1} - hV^0_{k}\right) .
 \end{split} \notag
\end{align}
In the above we used the abbreviation $D^{\pm}_k : = D_{\xi_k,\check{\mu}_k}^{\pm} $, where the maps  $D_{\xi, \mu}^{\pm}:\mathfrak{g}\to\mathfrak{g}^*$ are defined by the relation
\begin{equation}
\left.\frac{d}{d\varepsilon} \right|_{\varepsilon = 0} \left<\left(d\tau^{-1}_{\pm h(\xi + \varepsilon \rho)}\right)^*\mu, a\right>_{\mathfrak{g}^*\times \mathfrak{g}} =\left< D_{\xi,\mu}^{\pm}a, \rho\right>_{\mathfrak{g}^* \times \mathfrak{g}}, \quad \mbox{for all } \rho \in \mathfrak{g}
\,.
\end{equation}
Moreover we defined $\mathcal{A}_k: G \times \mathfrak{g}^* \to \mathfrak{g}^*$ as
\begin{equation}
\left<\mathcal{A}_k(g,V), \eta\right>  := \left<\Phi_k(g), \eta\right>  + \left.\frac{d}{d\varepsilon} \right|_{\varepsilon = 0}  \left< \Phi_k(g\exp\left(\varepsilon\eta\right)), V\right>\,, \notag
\end{equation}
for all $\eta\in\mathfrak{g}$. Then the gradient of the functional $\mathcal{J}_d$ defined in \eqref{disc_implicit_func} is given by
\begin{align}
\nabla_{\mu_0} \mathcal{J}_d &= -V^1_0, \qquad \nabla_{\nu_0}\mathcal{J}_d = - V^0_0. 
\end{align}

\subsubsection{Interpolation on $S^2$} 
We present simulations for $G=SO(3)$ and $V=\mathbb{R}^3$, where the group action is by matrix multiplication of vectors.
In this setting the $\diamond$ operator of \eqref{diamond_defn} becomes the standard cross product for $I,\omega\in\mathbb{R}^3$,
\begin{equation}\nonumber
I\diamond\omega = I\times\omega\in\mathbb{R}^3\cong\mathfrak{so}(3)\,.
\end{equation}
For our purposes we take the starting point $T_0 = (1,0,0)$, initial velocity $\xi_0^0 = \frac{5\pi}{2}(0, 0, 1)$, the node times $t_i = \frac{1}{5} i$ for $i = 1, \ldots ,5$ and target points 
\begin{equation} \notag
I_{t_1} = \left(
  \begin{array}{c}
    0\\
    1\\
    0\\
  \end{array}\right), \ I_{t_2} = \left(
  \begin{array}{c}
    0\\
    0\\
    1\\
  \end{array}\right), \ I_{t_3} =\frac{1}{\sqrt{2}}\left(
  \begin{array}{c}
    1\\
    0\\
    1\\
  \end{array}\right), \ I_{t_4} =  \frac{1}{\sqrt{2}}\left(
  \begin{array}{c}
    1\\
    1\\
    0\\
  \end{array}\right), \ I_{t_5} = \frac{1}{\sqrt{3}} \left(
  \begin{array}{c}
    1\\
    1\\
    1\\
  \end{array}\right).
\end{equation}
Note that all target points were chosen to have unit length so that we effectively perform interpolation on $S^2$. We refer to Figure \ref{Figure_TM_sphere_2} for an example interpolant and Figure \ref{mom_behaviour} for an illustration of the momentum behaviour.
\begin{SCfigure}
\centering
\includegraphics[scale=0.229]{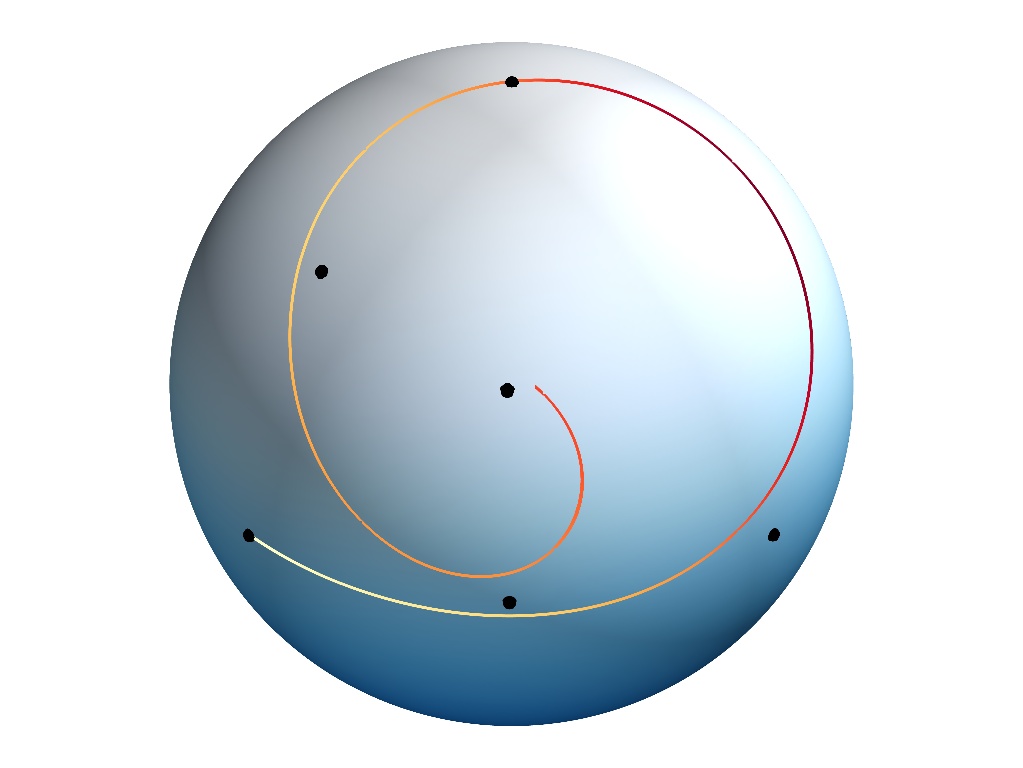}
\caption{\footnotesize Discrete-time trajectory planning on the sphere. The curve shown is a minimiser of the cost functional \eqref{Discrete_TM_functional} obtained numerically by following the procedures outlined in this section. Namely, the cost functional was restricted to solutions of the shooting equations \eqref{Discrete_HOHP_TM_equations}--\eqref{mu_discrete_node}, thus reducing the problem to the space of initial momenta, where a gradient descent was performed. The initial and target points were chosen as given above with tolerance parameter $\sigma = 0.025$. The colours represent the local speed along the curve in $SO(3)$, that is $\|\xi_k\|$ (red is large, white is small). }
\label{Figure_TM_sphere_2}
\end{SCfigure}
\begin{figure}[htb]
\begin{center}
\includegraphics[scale=0.55]{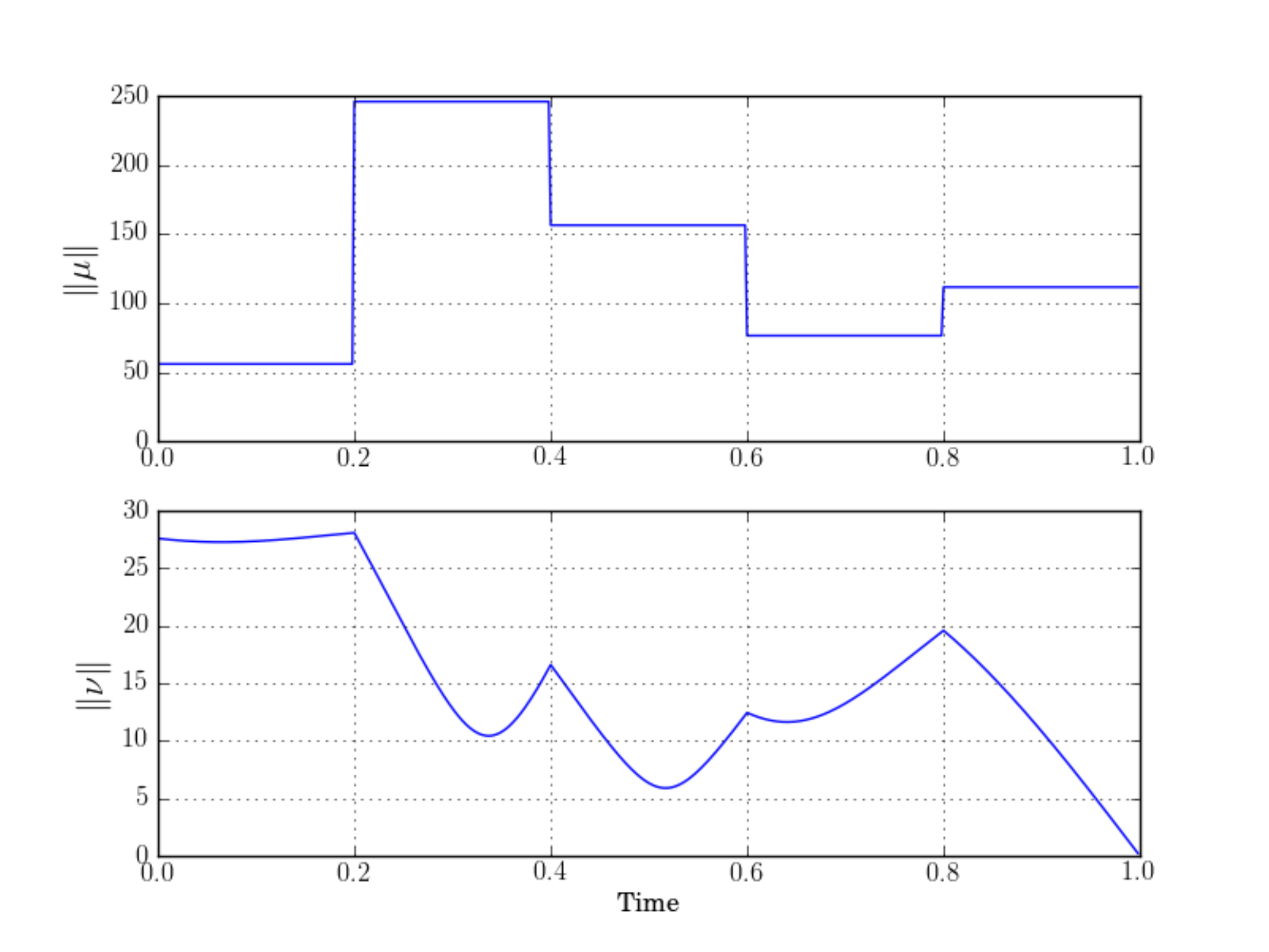}
\end{center}
\caption{\footnotesize Momentum norms. For the interpolating discrete cubic of Figure \ref{Figure_TM_sphere_2}, the plot shows the norms of the momenta $\mu_k$ and $\nu_k$. The norm of $\mu_k$ displays momentum kicks at node indices and exact conservation by discrete coadjoint motion for interior indices, as in equations \eqref{mu_discrete_interior} and \eqref{mu_discrete_node}. The norm of $\nu_k$ demonstrates continuity as found in the last equation of \eqref{Discrete_HOHP_TM_equations}. Both graphs respect terminal conditions \eqref{nu_discrete_endpoint} and \eqref{mu_discrete_endpoint}.}
\label{mom_behaviour}
\end{figure}
\section{Summary and Outlook}\label{Sum-sec}
This paper has discussed a structure-preserving numerical integration scheme for a class of higher-order mechanical systems. The structure-preserving properties arose as a a consequence of the variational nature of the discretisation method. 

Our starting point was a continuous time higher-order Hamilton--Pontryagin variational principle, where second-order kinematic constraints were enforced via Lagrange multipliers. These kinematic constraints were transferred to the discrete-time setting by combining Runge--Kutta and Runge--Kutta--Munthe-Kaas methods, in a similar way to the first-order treatment in \cite{BRMa09}. The resulting discrete flow maps were shown to preserve the symplectic form of higher-order mechanics and to respect momentum conservation in the presence of group symmetries. We proceeded with an application to a trajectory-planning problem familiar from inexact template matching in computational anatomy. The continous-time problem was presented in the higher-order Hamilton--Pontryagin framework, which led to a particularly clear exposition of the geometric properties of solution curves. These properties were finally shown to be inherited by the higher-order Hamilton--Pontryagin integration scheme. 

The treatment in the present paper was focused on mechanics based on Lagrangians depending on velocity and acceleration. It is however clear that a generalization to third, and higher, order follows in a straightforward manner. Other directions for further research include a rigorous analysis of the level of accuracy of the higher-order Hamilton--Pontryagin integration schemes. A related goal would be to include Runge--Kutta--Munthe-Kaas methods with a higher order of accuracy by choosing a larger truncation index in the sense of Theorem \ref{Thm_HLW}.

The type of trajectory-planning problem treated here is relevant in a number of finite-dimensional situations, whenever a path $g(t) T_0$ generated by Lie group transformations $g(t)$ is required to pass near prescribed target points $I_{t_i}$ at given times $t_i$.  We have treated the case of rotations, in which the Lie group is $SO(3)$. However, one may envision applications to many other problems involving different Lie groups. For example, one could imagine applying these methods in designing particle-beam optics, in which the rotations would be replaced by symplectic transformations acting on phase-space moments of the beam distribution function, such as its emittance \cite{Dragt2011}. 

\paragraph{Acknowledgements.} We thank M. Bruveris, C. J. Cotter, D. Mart\'in de Diego, L. Noakes, C. Tronci and F.-X. Vialard for encouraging comments and insightful remarks during the course of this work. 
This work was partially supported by the Royal Society of London Wolfson scheme 
and the European Research Council Advanced scheme.

{\footnotesize \bibliographystyle{alpha}
\bibliography{HOHP}}

\end{document}